\documentclass[11pt]{article}

\usepackage[margin=1.2in]{geometry}
\usepackage{amsmath,amssymb,amsthm,mathtools}
\usepackage{microtype,needspace}
\usepackage[colorlinks=true,linkcolor=blue,
citecolor=red,urlcolor=blue]{hyperref}
\usepackage[nameinlink,capitalize,noabbrev]{cleveref}

\newtheorem{theorem}{Theorem}[section]
\newtheorem{lemma}[theorem]{Lemma}
\newtheorem{proposition}[theorem]{Proposition}
\newtheorem{corollary}[theorem]{Corollary}
\numberwithin{equation}{section}

\newcommand{\R}{\mathbb R}
\newcommand{\D}{D^{1,2}(\R^N)}
\newcommand{\Cyl}{\mathcal C}
\newcommand{\SN}{S^{N-1}}
\newcommand{\Lrad}{\mathcal L^{\mathrm{rad}}_\nu}
\newcommand{\spn}{\operatorname{span}}

\title{Uniqueness and nondegeneracy for a critical Choquard equation with Hardy potential}

\author{Guangze Gu and Zhipeng Yang\thanks{Corresponding author: yangzhipeng326@163.com}}

\date{}

\AtEndDocument{%
  \par
  \bigskip
  \bigskip

 \noindent
  \textbf{Guangze Gu}\\[0.2em]
  \textsc{Department of Mathematics, Yunnan Normal University, Kunming 650500, China}\\
  \textsc{Yunnan Key Laboratory of Modern Analytical Mathematics and Applications, Kunming 650500, China}\\[0.3em]
  \textit{E-mail address}: \texttt{guangzegu@163.com}
  
\bigskip

 \noindent
  \textbf{Zhipeng Yang}\\[0.2em]
  \textsc{Department of Mathematics, Yunnan Normal University, Kunming 650500, China}\\
  \textsc{Yunnan Key Laboratory of Modern Analytical Mathematics and Applications, Kunming 650500, China}\\[0.3em]
  \textit{E-mail address}: \texttt{yangzhipeng326@163.com}
}

\begin{document}
\maketitle

\begin{abstract}
We study the critical Choquard equation
\[
-\Delta u-\mu |x|^{-2}u
=\bigl(|\cdot|^{2-N}*|u|^p\bigr)|u|^{p-2}u
\qquad\text{in }\R^N,
\]
where $|\cdot|^{2-N}$ is the Newton kernel, $p=\frac{N+2}{N-2}$, and $0<\mu<\frac{(N-2)^2}{4}$. For every $N\ge3$, we prove that positive solutions in $D^{1,2}(\R^N)$ are radially symmetric about the origin, unique up to dilation, and nondegenerate. Our proof uses normalized logarithmic derivatives to obtain a cooperative system. The resulting monotonicity principle underlies both the uniqueness and nondegeneracy arguments, and remains effective in the subquadratic range $1<p<2$, corresponding to $N\ge7$. Consequently, the linearized quadratic form has Morse index one and the sharp inequality admits a global quadratic remainder estimate. The best constant is strictly concave in the Hardy parameter on $\bigl[0,\frac{(N-2)^2}{4}\bigr)$, is of class $C^2$ with strictly negative second derivative in the interior, and satisfies two-sided bounds that agree to first order at $\mu=0$. 
\end{abstract}



\section{Introduction and main results}

We study the critical Choquard equation
\begin{equation}\label{eq1.1}
-\Delta u-\mu\frac{u}{|x|^2}
=\left(\int_{\R^N}\frac{|u(y)|^p}{|x-y|^{N-2}}\,dy\right)|u|^{p-2}u
\quad\text{in }\R^N,
\end{equation}
where
\[
N\ge3,\qquad p=\frac{N+2}{N-2},\qquad
0<\mu<\mu_H=\frac{(N-2)^2}{4}.
\]
We write the Newton kernel as $K_N(x)=|x|^{2-N}$. Since the nonlocal interaction and the Hardy term scale like the Dirichlet energy, the equation is invariant under dilations, whereas the Hardy potential breaks translation invariance. We classify positive solutions in $D^{1,2}(\R^N)$, determine the kernel of the linearized operator in this space, and prove a global quadratic remainder estimate for the sharp inequality. We also establish strict concavity and endpoint asymptotics of the best constant.
The proof combines normalized logarithmic derivatives with an auxiliary quotient to obtain a cooperative system. The vanishing of the Hamiltonian reduces the centered initial data to one parameter, and comparison of the associated trajectories proves uniqueness. The same comparison argument, applied to the linearized system, shows that the even radial kernel is trivial. A comparison between the solution and its Newton potential controls the singular factors when $1<p<2$, so the argument also covers dimensions $N\ge7$.

Let $\D$ be the completion of $C_c^\infty(\R^N)$ under the Dirichlet norm. Put
\[
a=\frac{N-2}{2},\qquad \nu=\sqrt{a^2-\mu},\qquad 2^*=\frac{2N}{N-2},
\]
and define
\[
\begin{split}
Q_\mu(u)&=\int_{\R^N}|\nabla u|^2\,dx-\mu\int_{\R^N}|x|^{-2}u^2\,dx,\\
D(u)&=\iint_{\R^N\times\R^N}\frac{|u(x)|^p|u(y)|^p}{|x-y|^{N-2}}\,dx\,dy.
\end{split}
\]
We use the same definition of $Q_\mu$ at $\mu=0$, with the Hardy term omitted. For $0\le\mu<\mu_H$, the best constant is
\begin{equation}\label{eq1.2}
S_\mu=\inf_{u\in\D\setminus\{0\}}\frac{Q_\mu(u)}{D(u)^{\frac{1}{p}}}.
\end{equation}
The Hardy, Sobolev, and Hardy-Littlewood-Sobolev inequalities imply
\[
Q_\mu(u)\ge\left(1-\frac{\mu}{\mu_H}\right)\|\nabla u\|_2^2,
\qquad D(u)^{\frac{1}{p}}\le C\|u\|_{2^*}^2,
\]
so $S_\mu>0$ \cite{LiebHLS,LiebLoss}. The equation and both integrals are invariant under
\[
u_\lambda(x)=\lambda^a u(\lambda x),\qquad\lambda>0.
\]
The quotient is also invariant under multiplication by a nonzero constant. We call a positive extremal normalized when its amplitude is chosen so that it solves \eqref{eq1.1}. A positive finite energy solution is a positive element of $\D$ satisfying \eqref{eq1.1} weakly.

The classical Choquard equation contains a positive mass term. In dimension three, Lieb \cite{Lieb77} proved existence and uniqueness, up to translation, of the minimizing solution of
\[
-\Delta Q+Q=(|\cdot|^{-1}*Q^2)Q.
\]
Lions \cite{LionsChoquard} developed the variational theory, and Ma and Zhao \cite{MaZhao} classified positive solitary solutions. Lenzmann \cite{Lenzmann} proved nondegeneracy of the nonrelativistic Hartree ground state as part of his study of the pseudorelativistic equation. For a general Riesz potential $I_\alpha(x)=c_{N,\alpha}|x|^{\alpha-N}$, with $0<\alpha<N$ and the usual positive normalization, the corresponding equation is
\[
-\Delta Q+Q=(I_\alpha*|Q|^q)|Q|^{q-2}Q.
\]
Moroz and Van Schaftingen \cite{MorozVS13} established existence, qualitative properties, and decay of ground states in the range
\[
\frac{N+\alpha}{N}<q<\frac{N+\alpha}{N-2}\qquad(N\ge3),
\]
and treated more general nonlinearities in \cite{MorozVS15}. Their survey \cite{MorozVSGuide} gives a broader account. For equations with a positive mass term, Xiang \cite{Xiang}, Seok \cite{Seok}, and Li \cite{Li22} proved uniqueness and nondegeneracy near particular powers or Riesz orders by perturbation from limiting equations.
Gu and Yang \cite{GuYang2026} proved uniqueness and nondegeneracy of
ground states for the mass-critical Hartree equation
\[
-\Delta Q+Q=(|\cdot|^{-2}*|Q|^2)Q
\quad\text{in }\R^N,\qquad N\ge3.
\]
Their uniqueness argument combines radial nondegeneracy with
continuation in a real dimension parameter. They also established
a global quadratic stability estimate for the sharp Hartree
interpolation inequality after normalizing the mass and gradient norm.

At the upper critical exponent, the free equation has explicit extremals. Put $q_\alpha=\frac{N+\alpha}{N-2}$. Since $\frac{2Nq_\alpha}{N+\alpha}=2^*$, the sharp Hardy-Littlewood-Sobolev and Sobolev inequalities apply with matching exponents. The bubbles of Aubin \cite{Aubin} and Talenti \cite{Talenti} attain equality in both inequalities \cite{LiebHLS}. Thus the positive extremals of the free critical quotient are
\[
Z_{c,\lambda,\xi}(x)
=c\left(\frac{\lambda}{\lambda^2+|x-\xi|^2}\right)^a,
\qquad c,\lambda>0,\quad\xi\in\R^N.
\]
After the amplitude is fixed, they solve
\[
-\Delta Z=(I_\alpha*Z^{q_\alpha})Z^{q_\alpha-1}.
\]
The integral method of Chen, Li, and Ou \cite{ChenLiOu} was the basis of later classifications of positive solutions to critical nonlocal equations, including the Choquard equation \cite{DuYang,GuoHuPengShuai}.

The corresponding nondegeneracy problem is to determine whether translations and dilation generate the linearized kernel. Let $L_{0,\alpha}$ denote the linearization at a normalized bubble $Z$. Du and Yang \cite{DuYang} and Giacomoni, Wei, and Yang \cite{GiacomoniWeiYang} proved nondegeneracy in partial parameter ranges. Li, Liu, Tang, and Xu \cite{LiLiuTangXu} treated the whole range $0<\alpha<N$ for bounded kernel elements vanishing at infinity. In that class,
\[
\left\{\phi\in L^\infty(\R^N):L_{0,\alpha}\phi=0,\ \lim_{|x|\to\infty}\phi(x)=0\right\}
=\spn\{aZ+x\cdot\nabla Z,\partial_{x_1}Z,\ldots,\partial_{x_N}Z\}.
\]
For $q_\alpha\ge2$, the same kernel characterization in $\D$ is proved in \cite[Proposition 1.4]{LiLiuTangXuThreshold} as part of the analysis of threshold dynamics. For the Newton kernel, this range includes $3\le N\le6$. With the Hardy potential, we show that the dilation derivative spans the kernel of the linearized operator in $\D$, whose elements may be singular at the origin.

For the local Hardy equation, Terracini \cite{Terracini} classified positive finite energy solutions of
\[
-\Delta w-\mu|x|^{-2}w=w^{\frac{N+2}{N-2}},\qquad 0<\mu<\mu_H.
\]
Up to dilation, the solution is
\[
w_\mu(r)=\frac{c_\mu}{r^{a-\nu}(1+r^{\frac{2\nu}{a}})^a},
\]
where $c_\mu>0$ is fixed by the equation. Dancer, Gladiali, and Grossi \cite{DancerGladialiGrossi} determined the linearized kernel, which for $0<\mu<\mu_H$ is generated by the dilation derivative. Robert \cite{Robert} established nondegeneracy for a broader class of Hardy-Sobolev equations. The local best constant is also explicit:
\[
S^{\mathrm{loc}}_\mu
=\inf_{u\in\D\setminus\{0\}}\frac{Q_\mu(u)}{\|u\|_{2^*}^2}
=S\left(1-\frac{\mu}{\mu_H}\right)^{\frac{N-1}{N}}.
\]
Here $S$ is the sharp Sobolev constant. In \eqref{eq1.1}, the solution is coupled to its Newton potential, and we classify the pair $(U,K_N*U^p)$.

For the critical Riesz equation with a Hardy potential,
\[
-\Delta u-\mu|x|^{-2}u=(I_\alpha*u^{q_\alpha})u^{q_\alpha-1},
\]
Guo and Tang \cite{GuoTang} proved variational existence, radial symmetry, monotonicity, and two-sided estimates at the origin and infinity under
\[
(N-4)_+<\alpha<N,\qquad 0<\mu<\mu_H.
\]
They obtained
\[
u(r)\asymp r^{-(a-\nu)}\quad(r\to0),\qquad
u(r)\asymp r^{-(a+\nu)}\quad(r\to\infty),
\]
where $\asymp$ denotes comparison by positive constants. The condition on $\alpha$ is equivalent to $q_\alpha>2$ and, for the Newton kernel ($\alpha=2$), covers exactly dimensions three, four, and five. Cao, Dai, and Li \cite{CaoDaiLi} proved symmetry and sharp asymptotic estimates for weighted quasilinear nonlocal equations. Their Theorem 1.7 includes the Newton equation in dimension six: in the case of the Laplacian (their $p=2$), the kernel is $|x|^{-4}$. Their semilinear results require the kernel exponent to be at most four, which for the Newton kernel means $N\le6$.

Guo, Gui, and Tang \cite{GuoGuiTang} studied existence with additional potentials and nonlinear perturbations. Ground states of related Choquard equations with Hardy potentials and critical nonlinearities were obtained in \cite{GuoTangZhang}. Li, Yang, and Zhou \cite{LiYangZhou} considered a coupled Hartree system with Hardy potentials. Solutions on $\R^N$ and sharp constants also determine concentration energies in critical problems on bounded domains. For Choquard equations, this connection is developed in \cite{GaoDaSilvaYangZhou,GaoYang} through the concentration-compactness method of Lions \cite{LionsCC}. For the Newton kernel in dimensions three, four, and five, Gu and Jevnikar \cite{GuJevnikar} proved
\[
S_0\left(1-\frac{\mu}{\mu_H}\right)^{\frac{N-1}{N}}<S_\mu<S_0.
\]
Theorem~\ref{Thm1.4} below extends these bounds to all $N\ge3$ and determines the rate at which $S_\mu$ vanishes at the Hardy threshold.

The exponent $p$ decreases with $N$:
\[
p>2\quad(3\le N\le5),\qquad p=2\quad(N=6),\qquad 1<p<2\quad(N\ge7).
\]
When $p<2$, the moving plane argument and the linearization involve the singular factor $U^{p-2}$. We estimate its product with the Newton potential. Write $\omega_N=|S^{N-1}|$ and $\kappa_N=(N-2)\omega_N$. By Lemma~\ref{Lem2.1},
\[
0<K_N*U^p\le\sqrt{\kappa_N}U,\qquad
(K_N*U^p)U^{p-2}\le\sqrt{\kappa_N}U^{p-1}\in L^{\frac{N}{2}}(\R^N)
\]
for every positive finite energy solution. The proof tests the equation for $\sqrt{\kappa_N}U-K_N*U^p$ against its negative part, so the bound is available at the start of the argument, before symmetry and decay are proved. We use this $L^{\frac{N}{2}}$ bound for the nonlinear differences in the moving plane argument and for the linearized potential.

At a positive solution $U$, write $\Lambda U=aU+x\cdot\nabla U$ and define the linearization by
\begin{equation}\label{eq1.3}
\begin{split}
L_\mu\phi={}&-\Delta\phi-\mu|x|^{-2}\phi
-(p-1)(K_N*U^p)U^{p-2}\phi\\
&-p\bigl(K_N*(U^{p-1}\phi)\bigr)U^{p-1}.
\end{split}
\end{equation}
Lemma~\ref{Lem4.1} proves Fr\'echet differentiability of the nonlinearity at $U$. We take the kernel of $L_\mu$ in $\D$, with the equation understood weakly.

\Needspace{12\baselineskip}
\begin{theorem}\label{Thm1.1}
Let $N\ge3$ and $0<\mu<\mu_H$. There is a normalized positive extremal $U_\mu$ of \eqref{eq1.2}. Every positive finite energy solution of \eqref{eq1.1} is of the form
\[
U_{\mu,\lambda}(x)=\lambda^aU_\mu(\lambda x),\qquad\lambda>0.
\]
Every nonzero real extremal of \eqref{eq1.2} is $cU_{\mu,\lambda}$ for some $c\in\R\setminus\{0\}$ and $\lambda>0$. Moreover,
\begin{equation}\label{eq1.4}
\ker_{\D}L_\mu=\spn\{\Lambda U_\mu\}.
\end{equation}
\end{theorem}

In particular, every positive finite energy solution is an extremal. Testing the equation against any nonzero solution gives $Q_\mu(U)=D(U)$, and \eqref{eq1.2} implies $Q_\mu(U)\ge S_\mu^{\frac{p}{p-1}}$. Thus every positive solution has the least energy
\[
\frac12Q_\mu(U)-\frac1{2p}D(U)
=\frac{p-1}{2p}S_\mu^{\frac{p}{p-1}}
\]
among nonzero finite energy solutions.

We fix the dilation by choosing $U_\mu$ so that
\[
v_\nu(t)=e^{at}U_\mu(e^t),\qquad\nu=\sqrt{a^2-\mu},
\]
is even. Lemma~\ref{Lem2.5} shows that this choice is unique. The centered solution is invariant under Kelvin inversion, and its leading terms at zero and infinity are given in Corollary~\ref{Cor3.5}.

For the quantitative statement, use the energy inner product
\[
\langle f,g\rangle_\mu=\int_{\R^N}\bigl(\nabla f\cdot\nabla g-\mu|x|^{-2}fg\bigr),
\qquad \|f\|_\mu^2=Q_\mu(f).
\]
Let $\mathfrak q_\mu$ be the quadratic form of \eqref{eq1.3} at $U_\mu$, and set
\[
\mathcal M_\mu=\{cU_{\mu,\lambda}:c\in\R,\ \lambda>0\},
\qquad d_\mu(u)=\inf_{z\in\mathcal M_\mu}\|u-z\|_\mu.
\]
Bianchi and Egnell \cite{BianchiEgnell} bounded the Sobolev deficit from below by the squared distance to the extremal family. Analogous estimates are known for the Caffarelli-Kohn-Nirenberg inequalities \cite{WeiWu} and for the Hardy-Sobolev inequality with a Hardy potential \cite{Chakraborty}. For the present inequality, we first derive a local estimate of this type from nondegeneracy and the second variation on the constraint $D(u)=D(U_\mu)$, and then extend the lower bound to all of $\D$ by concentration compactness.

\begin{theorem}\label{Thm1.2}
For each $N\ge3$ and $0<\mu<\mu_H$, the form $\mathfrak q_\mu$ has Morse index one. There are constants $c_\mu,\gamma_\mu,C_\mu>0$ and $\epsilon_\mu\in(0,1)$ such that
\begin{equation}\label{eq1.5}
\mathfrak q_\mu[\phi]\ge\gamma_\mu\|\phi\|_\mu^2
\quad\text{if }\langle\phi,U_\mu\rangle_\mu
=\langle\phi,\Lambda U_\mu\rangle_\mu=0,
\end{equation}
and
\begin{equation}\label{eq1.6}
c_\mu d_\mu(u)^2\le Q_\mu(u)-S_\mu D(u)^{\frac{1}{p}}
\le C_\mu d_\mu(u)^2
\end{equation}
whenever $u\in\D$ and $d_\mu(u)\le\epsilon_\mu\|u\|_\mu$.
\end{theorem}

The Morse index is the maximal dimension of a subspace on which the form is negative definite. In \eqref{eq1.5}, the first orthogonality condition defines the tangent space to $D(u)=D(U_\mu)$ at $U_\mu$, and the second removes the dilation direction. The constants may depend on $N$ and $\mu$.

\begin{corollary}\label{Cor1.3}
For each $N\ge3$ and $0<\mu<\mu_H$, there is $\hat c_\mu>0$ such that
\begin{equation}\label{eq1.7}
\hat c_\mu\, d_\mu(u)^2\le Q_\mu(u)-S_\mu D(u)^{\frac{1}{p}}\le d_\mu(u)^2
\qquad\text{for every }u\in\D .
\end{equation}
\end{corollary}

We next describe the dependence of $S_\mu$ on $\mu$ and its asymptotics at both endpoints.

\begin{theorem}\label{Thm1.4}
For every $N\ge3$, the function $\mu\mapsto S_\mu$ is strictly decreasing and strictly concave on $[0,\mu_H)$. It belongs to $C^2((0,\mu_H))$ and satisfies $S_\mu''<0$. For $0<\mu<\mu_H$,
\begin{equation}\label{eq1.8}
S_0\left(1-\frac{\mu}{\mu_H}\right)^{\frac{N-1}{N}}
<S_\mu<
S_0\left(1-\frac{N-1}{N}\frac{\mu}{\mu_H}\right).
\end{equation}
Consequently,
\begin{equation}\label{eq1.9}
S_\mu=S_0\left(1-\frac{4(N-1)}{N(N-2)^2}\mu\right)+O(\mu^2)
\quad\text{as }\mu\to0^+.
\end{equation}
At the Hardy threshold,
\begin{equation}\label{eq1.10}
\begin{split}
S_\mu&=C_N(\mu_H-\mu)^{\frac{N}{N+2}}
\bigl[1+b_N(\mu_H-\mu)+o(\mu_H-\mu)\bigr],\\
b_N&=\frac{8(N+2)}{(N-2)^3(N+4)},
\end{split}
\end{equation}
where
\[
C_N=\omega_N^{1-\frac{2}{p}}\left(\frac a2\right)^{\frac{1}{p}}I_p^{\frac{p-1}{p}},
\qquad
I_p=\frac{p^{\frac{p}{p-1}}\sqrt\pi}{p-1}
\frac{\Gamma\!\left(\frac{p}{p-1}\right)}{\Gamma\!\left(\frac{p}{p-1}+\frac{1}{2}\right)}.
\]
\end{theorem}

Since $\frac{N}{N+2}<\frac{N-1}{N}$, the constant $S_\mu$ tends to zero more slowly than the local Hardy-Sobolev constant. We obtain the exponent and the correction in \eqref{eq1.10} by passing to the limit in the centered profiles. With $M_N=\frac{2\omega_N}{a}$, Proposition~\ref{Prop5.2} proves
\[
\left(\frac{M_N}{\nu^2}\right)^{\frac{1}{2p-2}}v_\nu\!\left(\frac{\cdot}{\nu}\right)
\longrightarrow w_*\quad\text{strongly in }H^2(\R)
\]
as $\nu\to0^+$, where
\[
w_*(s)=\frac{p^{\frac{1}{2p-2}}}{\cosh^{\frac{1}{p-1}}((p-1)s)},
\qquad -w_*''+w_*=w_*^{2p-1}.
\]
The proof of strict concavity uses the following parameter dependence of the centered solution family. Its derivative determines the second derivative of $S_\mu$.

\begin{corollary}\label{Cor1.5}
The centered normalized profile is continuously differentiable in $H^2(\R)$ with respect to $\mu\in(0,\mu_H)$. The corresponding family $U_\mu$ is continuously differentiable in $\D$, and
\begin{equation}\label{eq1.11}
S_\mu'=-\frac{\displaystyle\int_{\R^N}|x|^{-2}U_\mu^2\,dx}{D(U_\mu)^{\frac{1}{p}}}.
\end{equation}
\end{corollary}

We outline the proof. For $V=K_N*U^p$, subtracting the two weak equations gives
\[
-\Delta(\sqrt{\kappa_N}U-V)
+\sqrt{\kappa_N}U^{p-1}(\sqrt{\kappa_N}U-V)
=\sqrt{\kappa_N}\mu|x|^{-2}U.
\]
Testing with the negative part proves $V\le\sqrt{\kappa_N}U$. The same bound for reflected functions gives the difference estimates needed for radial symmetry. When $\mu=0$, testing with $\sqrt{\kappa_N}U-V$ instead gives equality and reduces the Newton system to a scalar critical equation.

For radial solutions, the change of variables $U(r)=r^{-a}v(t)$, $V(r)=r^{-a}W(t)$, and $t=\log r$ gives
\begin{equation}\label{eq1.12}
-v''+\nu^2v=Wv^{p-1},\qquad
-W''+a^2W=\kappa_Nv^p.
\end{equation}
Here $W=k*v^p$, with $k(t)=\omega_Ne^{-a|t|}$. A second moving plane argument establishes a common symmetry center on the line. In dimension six, $p=2$ and the substitution $f=\sqrt{2\kappa_6}v$, $g=W$ transforms \eqref{eq1.12} into
\[
-f''+\nu^2f=fg,\qquad -g''+a^2g=\frac12f^2,
\]
the system studied by Lopes \cite{Lopes}. His uniqueness proof proceeds by continuation, using nondegeneracy on the even subspace and uniform a priori estimates. He also recalled that the linearized kernel is generated by the translation derivative. Whereas scalar uniqueness arguments such as those in \cite{Kwong} concern a single equation, the normalized variables below give a direct comparison for the coupled system \eqref{eq1.12} and its linearization in every dimension $N\ge3$.

The main ingredient in the radial analysis is a normalization of the logarithmic derivatives. For a centered profile, we introduce
\[
X=-\frac{v'}{v},\qquad Y=-\frac{W'}{W},\qquad
q=\frac{Wv^{p-2}}{\nu^2-X^2}.
\]
The Green representations give $|X|<\nu$ and $|Y|<a$, while the vanishing of the Hamiltonian yields $q>\frac p2$ and allows us to eliminate $\frac{\kappa_Nv^p}{W}$. After rescaling the independent variable and setting $x=\frac{X}{\nu}$ and $z=-\frac{Y}{a}$, we obtain the cooperative system \eqref{eq3.10}. Its centered initial data are determined by a single parameter $q_0$. Two solutions with different values of $q_0$ would remain strictly ordered, and the ratio of their $q$ components would be strictly increasing. Lemma~\ref{Lem3.2}, however, gives the same positive limit for both components. This proves uniqueness without an explicit formula for the profile.

The same comparison argument shows that the even radial kernel is trivial. A nonzero even solution of the linearized equations determines a nonzero variation $\delta q_0$ of the centered initial data. After choosing its sign, the variational system makes $\frac{\delta q}{q}$ positive and strictly increasing. Lemma~\ref{Lem4.3} proves directly from the linearized equations that $\frac{\delta q}{q}\to0$, giving a contradiction. The odd radial kernel is then identified as $\spn\{v_\nu'\}$ by a ground state representation on the half-line. A separate identity gives
\[
\mathfrak q_\ell[f]\ge\left(\frac{2\mu}{a+\nu}+\lambda_\ell-\lambda_1\right)\|f\|_{L^2(\R,dt)}^2,
\qquad \lambda_\ell=\ell(\ell+N-2),\quad\ell\ge1,
\]
for every nonradial spherical harmonic component. The constrained second variation then yields Theorem~\ref{Thm1.2}. Since $S_\mu<S_0$, minimizing sequences are compact up to dilation, which gives Corollary~\ref{Cor1.3}.

At the Hardy threshold, rescaling gives the convolution kernel $\rho_\nu(s)=\frac{a}{2\nu}e^{-\frac{a|s|}{\nu}}$, whose mass is one. Centered minimizers converge to the limiting profile $w_*$, and the Fourier multiplier of $\rho_\nu$ gives the correction of order $\nu^2$. Testing with the limiting profile gives an upper bound, and the convergent minimizers give a matching lower bound. Qualitative $H^2$ convergence of the profiles suffices for \eqref{eq1.10}.

The Newton kernel is used in an essential way. By Lemma~\ref{Lem2.3}, the interaction on the cylinder is $\kappa_N(-\partial_t^2-\Delta_\theta+a^2)^{-1}$, so radial solutions satisfy the second-order system \eqref{eq1.12}. Its conserved Hamiltonian and the comparison principle for the associated cooperative system are used in Sections~3 and~4. For general Riesz orders, the interaction does not give the second-order equation for the potential in \eqref{eq1.12}, and our argument does not apply directly.

Section~2 establishes comparison, symmetry, and existence. Section~3 proves uniqueness and the profile asymptotics. Section~4 treats the linearized kernel, Morse index, local and global stability estimates, and parameter differentiability. Section~5 proves the statements about the optimal constant and its endpoint expansion.

\section{Energy estimates, symmetry, and existence}

We write $\|\cdot\|_q$ for the $L^q$ norm on the indicated domain. The identities
\[
p\frac{2N}{N+2}=2^*=\frac{2N}{N-2},\qquad
(p-1)\frac N2=2^*
\]
imply that the power map $u\mapsto |u|^p$ is continuously differentiable from $L^{2^*}$ to $L^{\frac{2N}{N+2}}$. The Hardy-Littlewood-Sobolev inequality then gives $D\in C^1(\D)$, with
\[
D'(u)[h]=2p\int_{\R^N}(K_N*|u|^p)|u|^{p-2}uh.
\]
The nonlinear term in \eqref{eq1.1} belongs to $L^{\frac{2N}{N+2}}$, and the Hardy potential defines an element of the dual of $\D$. Since a point has zero $2$-capacity for $N\ge3$, $C_c^\infty(\R^N\setminus\{0\})$ is dense in $\D$. Thus the weak identity extends from tests supported away from the origin to all tests in $\D$.

\begin{lemma}\label{Lem2.1}
Let $U\in\D$ be a nonzero nonnegative weak solution of \eqref{eq1.1}, and let $V=K_N*U^p$. Then $U,V$ are positive and smooth on $\R^N\setminus\{0\}$, $V\in\D$, and
\begin{equation}\label{eq2.1}
0<V\le\sqrt{\kappa_N}U,
\qquad 0<VU^{p-2}\le\sqrt{\kappa_N}U^{p-1}\in L^{\frac{N}{2}}(\R^N).
\end{equation}
For $\mu=0$, \eqref{eq2.1} holds with $V=\sqrt{\kappa_N}U$.
\end{lemma}

\begin{proof}
Since $U^p\in L^{\frac{2N}{N+2}}$, the Newton potential belongs to $L^{2^*}\cap\D$, satisfies $-\Delta V=\kappa_NU^p$, and obeys
\[
\|\nabla V\|_2\le C\|U^p\|_{\frac{2N}{N+2}}.
\]
Here $C$ depends only on $N$. Put $h=\sqrt{\kappa_N}U-V$. Subtracting the two weak equations gives
\[
-\Delta h+\sqrt{\kappa_N}U^{p-1}h
=\sqrt{\kappa_N}\mu|x|^{-2}U.
\]
Since $U^{p-1}\in L^{\frac{N}{2}}$, the zeroth order term is a continuous bilinear form on $\D$. The Hardy inequality controls the right side. Testing with $h_-=\max\{-h,0\}$ and multiplying by $-1$ gives
\[
\int|\nabla h_-|^2+\sqrt{\kappa_N}\int U^{p-1}h_-^2
=-\sqrt{\kappa_N}\mu\int|x|^{-2}Uh_-\le0.
\]
It follows that $h\ge0$. If $\mu=0$, testing with $h$ yields $h=0$. The integral representation gives $V>0$ almost everywhere, so the comparison implies $U>0$ almost everywhere and the second estimate in \eqref{eq2.1}.

To prove local regularity away from the origin, write on an open set $\Omega\Subset\R^N\setminus\{0\}$
\[
-\Delta U=A(x)U,\qquad A(x)=\mu|x|^{-2}+VU^{p-2}\in L^{\frac{N}{2}}.
\]
The local iteration of Brezis and Kato \cite{BrezisKato} gives $U\in L^q_{\mathrm{loc}}$ for every finite $q$. 

The equation for $V$ gives $V\in W^{2,q}_{\mathrm{loc}}$ for every finite $q$. Since $V$ is locally bounded and $U^{p-1}$ belongs locally to every finite $L^q$ space, the first equation gives the same regularity for $U$. Both functions are continuous and, by the strong maximum principle, strictly positive. Elliptic regularity then gives smoothness away from the origin.
\end{proof}

\begin{proposition}\label{Prop2.2}
Every positive finite energy solution of \eqref{eq1.1}, with $0<\mu<\mu_H$, is radial about the origin and strictly decreasing as a function of the radius. The same assertions hold for its Newton potential.
\end{proposition}

\begin{proof}
Fix a direction, which we take to be the first coordinate direction. For $\lambda<0$, put
\[
\Sigma_\lambda=\{x_1<\lambda\},\qquad
x^\lambda=(2\lambda-x_1,x'),\qquad
U_\lambda(x)=U(x^\lambda),\quad V_\lambda(x)=V(x^\lambda).
\]
Let $\xi=(U-U_\lambda)_+$ and $\zeta=(V-V_\lambda)_+$ in $\Sigma_\lambda$, extended by zero outside this half-space, and set
\[
E_\lambda=\{x\in\Sigma_\lambda:U(x)>U_\lambda(x)\}.
\]
Both differences have zero trace on the reflecting plane, and truncation followed by zero extension preserves $\D$. Thus $\xi,\zeta\in\D$ are admissible in the weak difference equations on $\Sigma_\lambda$, once integrability at the reflected singularity is checked below.

On $E_\lambda$, decompose the nonlinear difference as
\[
VU^{p-1}-V_\lambda U_\lambda^{p-1}
=(V-V_\lambda)U^{p-1}
+V_\lambda(U^{p-1}-U_\lambda^{p-1}).
\]
For $1<p<2$, concavity gives
\[
U^{p-1}-U_\lambda^{p-1}
\le(p-1)U_\lambda^{p-2}(U-U_\lambda).
\]
By Lemma~\ref{Lem2.1}, $V_\lambda U_\lambda^{p-2}\le\sqrt{\kappa_N}U_\lambda^{p-1}\le\sqrt{\kappa_N}U^{p-1}$ on this set. For $p\ge2$, the mean value theorem gives the same conclusion by evaluating the derivative at the larger endpoint and using $V_\lambda\le\sqrt{\kappa_N}U_\lambda$. Thus, in both cases,
\begin{equation}\label{eq2.2}
VU^{p-1}-V_\lambda U_\lambda^{p-1}
\le U^{p-1}\zeta+\sqrt{\kappa_N}(p-1)U^{p-1}\xi
\quad\text{on }E_\lambda.
\end{equation}

For $x\in\Sigma_\lambda$, $|x^\lambda|<|x|$. Subtraction of the first equations gives
\[
\begin{split}
&-\Delta(U-U_\lambda)-\mu|x|^{-2}(U-U_\lambda)\\
&\quad=VU^{p-1}-V_\lambda U_\lambda^{p-1}
+\mu\bigl(|x|^{-2}-|x^\lambda|^{-2}\bigr)U_\lambda\\
&\quad\le VU^{p-1}-V_\lambda U_\lambda^{p-1}.
\end{split}
\]
Test the two difference equations with $\xi$ and $\zeta$, respectively. Using \eqref{eq2.2}, the Hardy inequality, and $U^p-U_\lambda^p\le0$ on $\Sigma_\lambda\setminus E_\lambda$, we obtain
\begin{equation}\label{eq2.3}
\begin{split}
\left(1-\frac\mu{\mu_H}\right)\|\nabla\xi\|_2^2
&\le C_{N,p}\int_{E_\lambda}U^{p-1}(\xi^2+\xi\zeta),\\
\|\nabla\zeta\|_2^2
&\le\kappa_Np\int_{E_\lambda}U^{p-1}\xi\zeta.
\end{split}
\end{equation}
The tests are admissible at the reflected singularity $x_\lambda=(2\lambda,0,\ldots,0)$. Indeed,
\[
\int|x-x_\lambda|^{-2}U_\lambda\xi
\le\left(\int|x-x_\lambda|^{-2}U_\lambda^2\right)^{\frac{1}{2}}
\left(\int|x-x_\lambda|^{-2}\xi^2\right)^{\frac{1}{2}}<\infty
\]
by the translated Hardy inequality. The weak identities follow by density using cutoffs at the singular points and infinity.

Applying the H\"older and Sobolev inequalities to \eqref{eq2.3}, we obtain
\begin{equation}\label{eq2.4}
\begin{split}
\|\nabla\xi\|_2^2+\|\nabla\zeta\|_2^2
\le C_{N,\mu}\|U^{p-1}\|_{L^{\frac{N}{2}}(E_\lambda)}
\bigl(\|\nabla\xi\|_2^2+\|\nabla\zeta\|_2^2\bigr).
\end{split}
\end{equation}
The constant in \eqref{eq2.4} is independent of $\lambda$. Since $U^{p-1}\in L^{\frac{N}{2}}$, its norm on $\Sigma_\lambda$ tends to zero as $\lambda\to-\infty$. The estimate then gives $U\le U_\lambda$ and $V\le V_\lambda$ in $\Sigma_\lambda$ for all sufficiently negative $\lambda$.

Let $\lambda_0\le0$ be the supremum of the values $\lambda\le0$ such that both inequalities hold in $\Sigma_\eta$ for every $\eta\le\lambda$. Suppose $\lambda_0<0$. By continuity, the inequalities hold at $\lambda_0$ away from the singular points. Thus $w=U_{\lambda_0}-U$ and $z=V_{\lambda_0}-V$ are nonnegative on $\Sigma_{\lambda_0}\setminus\{x_{\lambda_0}\}$. Moreover,
\[
\begin{split}
-\Delta w={}&\mu|x|^{-2}w
+\mu\bigl(|x^{\lambda_0}|^{-2}-|x|^{-2}\bigr)U_{\lambda_0}\\
&+V_{\lambda_0}U_{\lambda_0}^{p-1}-VU^{p-1}>0,
\qquad
-\Delta z=\kappa_N(U_{\lambda_0}^p-U^p)\ge0.
\end{split}
\]
The strong maximum principle first gives $w>0$ and then $z>0$ on this punctured half-space.

Choose a compact set $K\subset\Sigma_{\lambda_0}\setminus\{x_{\lambda_0}\}$ such that $\int_{\Sigma_{\lambda_0}\setminus K}U^{2^*}$ is sufficiently small. The strict inequalities on $K$ persist when $\lambda>\lambda_0$ is close to $\lambda_0$, so $E_\lambda\cap K=\varnothing$. Also, $\int_{\Sigma_\lambda\setminus\Sigma_{\lambda_0}}U^{2^*}\to0$ as $\lambda\downarrow\lambda_0$. The coefficient in \eqref{eq2.4} is therefore less than one for all such $\lambda$. Both comparisons continue past $\lambda_0$, a contradiction. Hence $\lambda_0=0$.

Moving the plane from the opposite direction gives symmetry across the plane through the origin. Since the direction is arbitrary, both functions are radial. The strict comparison at every $\lambda<0$ gives strict radial decrease by reflecting two points on a ray across their midpoint.
\end{proof}

Write $x=r\theta$ and $t=\log r$, and let $\Cyl=\R\times\SN$ have product measure $dt\,d\theta$. For $u\in\D$, define its cylinder representative by
\begin{equation}\label{eq2.5}
u(r\theta)=r^{-a}v(t,\theta).
\end{equation}
We use the convention $-\Delta_\theta\ge0$. Direct integration for functions supported away from the origin, followed by density, gives
\begin{equation}\label{eq2.6}
Q_\mu(u)=\int_\Cyl\bigl(|v_t|^2+|\nabla_\theta v|^2+\nu^2v^2\bigr)\,dt\,d\theta.
\end{equation}
With $a^2$ instead of $\nu^2$, this is the Dirichlet norm identity. In particular, \eqref{eq2.5} identifies $\D$ with $H^1(\Cyl)$. This cylinder formulation is classical; see, e.g., \cite{CatrinaWang}.

The identity
\[
|e^t\theta-e^s\eta|^2=2e^{t+s}\bigl(\cosh(t-s)-\theta\cdot\eta\bigr)
\]
and $ap=a+2$ give
\begin{equation}\label{eq2.7}
D(u)=\iint_{\Cyl\times\Cyl}\mathcal K(t-s,\theta,\eta)
|v(t,\theta)|^p|v(s,\eta)|^p\,dt\,d\theta\,ds\,d\eta,
\end{equation}
where $\mathcal K(\tau,\theta,\eta)=[2(\cosh\tau-\theta\cdot\eta)]^{-a}$. Denote the corresponding integral operator by $T$. Let $e_1=(1,0,\ldots,0)\in\SN$.

\begin{lemma}\label{Lem2.3}
For radial functions on the cylinder, the convolution kernel and its mass are
\begin{equation}\label{eq2.8}
k(t)=\int_\SN\mathcal K(t,e_1,\eta)\,d\eta
=\omega_Ne^{-a|t|},\qquad
M_N=\int_\R k=\frac{2\omega_N}{a}=\frac{\kappa_N}{a^2}.
\end{equation}
On $L^2(\Cyl)$,
\begin{equation}\label{eq2.9}
T=\kappa_N(-\partial_t^2-\Delta_\theta+a^2)^{-1},\qquad
\|T\|_{L^2\to L^2}\le M_N.
\end{equation}
In particular, $T$ is self-adjoint and maps $L^2(\Cyl)$ continuously into $H^2(\Cyl)$.
\end{lemma}

\begin{proof}
Newton's shell formula is
\[
\int_\SN|re_1-s\eta|^{2-N}\,d\eta
=\omega_N\max\{r,s\}^{2-N}.
\]
The distance identity gives \eqref{eq2.8}. Translation and rotation invariance make every row integral of $\mathcal K$ equal to $M_N$. By symmetry, the column integrals are the same, and Schur's test gives the $L^2$ bound. To obtain the resolvent formula in \eqref{eq2.9}, take $f\in C_c^\infty(\Cyl)$ and apply the Newton potential to $r^{-a-2}f(\log r,\theta)$. The change of variables in \eqref{eq2.7} identifies its cylinder representative with $Tf$. Since $-\Delta K_N=\kappa_N\delta_0$,
\[
(-\partial_t^2-\Delta_\theta+a^2)Tf=\kappa_N f.
\]
By coercivity, the solution is unique in $H^1(\Cyl)$. Fourier transformation in $t$ and spherical harmonic expansion yield the resolvent identity and the $H^2$ estimate. Density extends these identities to $L^2(\Cyl)$.
\end{proof}

For radial cylinder functions, set
\begin{equation}\label{eq2.10}
\begin{split}
A_\nu(v)&=\int_\R(v'^2+\nu^2v^2),\qquad
B(v)=\int_\R(k*|v|^p)|v|^p,\\
s_\nu&=\inf_{v\in H^1(\R)\setminus\{0\}}
\frac{A_\nu(v)}{B(v)^{\frac{1}{p}}}.
\end{split}
\end{equation}
A nonnegative function on the line will be called centered if it is even and nonincreasing on $(0,\infty)$.

\begin{proposition}\label{Prop2.4}
For $0<\mu<\mu_H$,
\begin{equation}\label{eq2.11}
S_\mu=\omega_N^{1-\frac{1}{p}}s_\nu,
\end{equation}
and both infima are attained. For every $0<\nu\le a$, the one-dimensional infimum is attained by a centered function. After amplitude normalization, such a minimizer is a positive $H^1$ solution of
\begin{equation}\label{eq2.12}
-v''+\nu^2v=(k*v^p)v^{p-1}\quad\text{on }\R.
\end{equation}
\end{proposition}

\begin{proof}
Passing from $u$ to $|u|$ preserves the quotient. Symmetric decreasing rearrangement decreases the Dirichlet integral and increases both the Hardy integral and $D$ \cite{LiebLoss}. Since $Q_\mu(u)>0$ for $u\ne0$, the infimum is unchanged by restriction to nonnegative radial functions. Equations \eqref{eq2.6} and \eqref{eq2.7} give \eqref{eq2.11}.

On the line, symmetric decreasing rearrangement decreases $A_\nu$ and increases $B$. Young's inequality gives $B(v)\le M_N\|v\|_{2p}^{2p}$. The embedding $H^1(\R)\subset L^{2p}(\R)$ and the equivalence of $A_\nu^{\frac{1}{2}}$ with the $H^1$ norm imply $s_\nu>0$. Choose a centered minimizing sequence with $B=1$. It is bounded in $H^1$. For a centered function $h$, one has
\begin{equation}\label{eq2.13}
2rh(r)^2\le\|h\|_2^2,\qquad
\int_{|t|>R}|h|^q
\le\left(\frac{\|h\|_2}{\sqrt{2R}}\right)^{q-2}\|h\|_2^2
\quad(q>2).
\end{equation}
Local compactness and \eqref{eq2.13}, with $q=2p$, give a weakly convergent subsequence in $H^1$ that converges strongly in $L^{2p}(\R)$. Its limit remains centered. The $p$th powers converge in $L^2$, so Young's inequality implies that the limit satisfies $B=1$. Weak lower semicontinuity proves attainment.

A minimizer with $B=1$ satisfies \eqref{eq2.12} with multiplier $s_\nu$ on the right. Multiplication by $s_\nu^{\frac{1}{2p-2}}$ gives the normalized equation. Its right side belongs to $L^2$, since $v^p$ and $k*v^p$ lie in $L^2$ and $v^{p-1}$ is bounded. The exponential Green representation gives strict positivity. This argument also applies at $\nu=a$.
\end{proof}

\begin{lemma}\label{Lem2.5}
Let $0<\nu\le a$ and let $v$ be a positive $H^1(\R)$ solution of \eqref{eq2.12}. Put $W=k*v^p$. Then $v\in H^3(\R)$, $W\in H^2(\R)$, and the pair solves
\begin{equation}\label{eq2.14}
-v''+\nu^2v=Wv^{p-1},\qquad
-W''+a^2W=\kappa_Nv^p.
\end{equation}
There is a unique common symmetry center. After translating this center to zero, both functions are even and strictly decreasing on $(0,\infty)$, and
\begin{equation}\label{eq2.15}
W\le\sqrt{\kappa_N}v,\qquad
|v'|<\nu v,\qquad |W'|<aW,\qquad v''(0)<0.
\end{equation}
Both functions and their first derivatives tend to zero at either end of the line.
\end{lemma}

\begin{proof}
Since $v\in H^1(\R)$, we have $v^p\in L^2$. Since $k(t)=\frac{\kappa_N}{2a}e^{-a|t|}$ is $\kappa_N$ times the Green function of $-\frac{d^2}{dt^2}+a^2$, the convolution $W$ satisfies the second equation and belongs to $H^2$. Also $Wv^{p-1}\in L^2$, so the first equation gives $v\in H^2$. Both functions and their first derivatives tend to zero at either end of the line. For $h=\sqrt{\kappa_N}v-W$, subtraction gives
\[
-h''+(a^2+\sqrt{\kappa_N}v^{p-1})h
=\sqrt{\kappa_N}(a^2-\nu^2)v\ge0.
\]
Testing with the negative part of $h$ yields $W\le\sqrt{\kappa_N}v$. Hence
\begin{equation}\label{eq2.16}
0<R:=Wv^{p-2}\le\sqrt{\kappa_N}v^{p-1},\qquad
R\in L^\infty(\R),\qquad R(t)\to0\quad(|t|\to\infty).
\end{equation}
Now $(Wv^{p-1})'=W'v^{p-1}+(p-1)Rv'\in L^2$, and the first equation gives $v\in H^3$.

To prove symmetry, reflect the line about $t=\lambda$ and put $v_\lambda(t)=v(2\lambda-t)$ and $W_\lambda(t)=W(2\lambda-t)$. On $(-\infty,\lambda)$, let $\xi=(v-v_\lambda)_+$, $\zeta=(W-W_\lambda)_+$, and $E_\lambda=\{t<\lambda:v(t)>v_\lambda(t)\}$. Extend $\xi$ and $\zeta$ by zero. Since the original and reflected pairs satisfy $W\le\sqrt{\kappa_N}v$, the estimate in \eqref{eq2.2} applies. Testing the difference equations gives
\begin{equation}\label{eq2.17}
\begin{split}
\mathcal A_\lambda
&:=\int_\R(\xi'^2+\nu^2\xi^2+\zeta'^2+a^2\zeta^2)\\
&\le C\int_{E_\lambda}v^{p-1}(\xi^2+\xi\zeta).
\end{split}
\end{equation}
For sufficiently negative $\lambda$, $\sup_{t<\lambda}v(t)$ is small enough to absorb the right side, so both comparisons hold. Let $\lambda_0$ be the supremum of the values $\lambda$ for which they hold at every plane $\eta\le\lambda$.

Choose $t_0$ so that $v^{p-1}$ is small outside $[-t_0,t_0]$. The contribution from this complement can be absorbed in \eqref{eq2.17}. For a measurable set $E\subset[-t_0,t_0]$, the one-dimensional Sobolev embedding gives
\[
\int_Ev^{p-1}(\xi^2+|\xi\zeta|)
\le C_{t_0}|E|\bigl(\|\xi\|_{H^1}^2+\|\zeta\|_{H^1}^2\bigr)
\le C_{t_0,a,\nu}|E|\mathcal A_\lambda.
\]
The contribution from a set of sufficiently small measure can therefore be absorbed. At a finite limiting plane, the differences $w=v_{\lambda_0}-v$ and $z=W_{\lambda_0}-W$ are nonnegative and satisfy
\[
-w''+\nu^2w\ge0,\qquad -z''+a^2z\ge0.
\]
By the strong maximum principle, either both differences are positive in the interior or both vanish. If one vanishes identically, \eqref{eq2.14} forces the other to vanish. Suppose both are positive at $\lambda_0$. Enlarge $t_0$ so that $t_0>|\lambda_0|+1$. For small $\delta>0$, they have a positive lower bound on $[-t_0,\lambda_0-\delta]$, which persists for nearby planes. Thus $E_\lambda\cap[-t_0,t_0]\subset(\lambda_0-\delta,\lambda)$ for $\lambda>\lambda_0$ close enough to $\lambda_0$. This interval can be made sufficiently short for the preceding estimate to apply. Combined with the estimate outside $[-t_0,t_0]$, it continues both comparisons past $\lambda_0$, a contradiction. A finite limiting plane is therefore a common symmetry plane.

If the plane could move to $+\infty$, then $v$ would be nondecreasing on $\R$, contrary to positivity and decay at both ends. Thus the limiting plane is finite. Two distinct common centers would make $v$ periodic, contradicting its decay. The reflection inequalities give monotonicity away from the unique center.

Translate the center to zero. By the layer cake formula, an even nonincreasing positive function is an integral of indicators of centered intervals. For $t,r>0$,
\[
\frac d{dt}\bigl(k*\mathbf1_{[-r,r]}\bigr)(t)=k(t+r)-k(t-r)<0.
\]
Applying this identity to $W=k*v^p$ gives $W'(t)<0$. Thus $F=Wv^{p-1}$ is strictly decreasing on $(0,\infty)$. The same calculation for $v=G_\nu*F$, where $G_\nu(t)=\frac{e^{-\nu|t|}}{2\nu}$, gives $v'(t)<0$.

In each Green representation, the integrals over the two half-lines separated by $t$ are strictly positive. Differentiation subtracts these terms and multiplies by $\nu$ or $a$. This gives $|v'|<\nu v$ and $|W'|<aW$. At the center,
\[
\nu^2v(0)=\frac\nu2\int_\R e^{-\nu|s|}F(s)\,ds<F(0),
\]
so the first equation gives $v''(0)<0$.
\end{proof}

\section{Uniqueness of positive solutions}

We now prove uniqueness for positive $H^1$ solutions of the cylinder system. The decay estimates rely on $W\le\sqrt{\kappa_N}v$ from Lemma~\ref{Lem2.5}, which holds in every dimension $N\ge3$.

Fix $0<\nu\le a$ and a positive $H^1(\R)^2$ solution $(v,W)$ of \eqref{eq2.14}. The second equation and coercivity give $W=k*v^p$. By Lemma~\ref{Lem2.5}, we may center the pair at zero. Its Hamiltonian satisfies
\begin{equation}\label{eq3.1}
H=\frac{v'^2}{2}+\frac{W'^2}{2\kappa_Np}
-\frac{\nu^2v^2}{2}-\frac{a^2W^2}{2\kappa_Np}
+\frac1pWv^p=0.
\end{equation}
Indeed, differentiation using \eqref{eq2.14} gives $H'=0$, and decay at infinity implies $H=0$.

\begin{lemma}\label{Lem3.1}
There are constants $A>0$ and $\eta>0$ such that, as $t\to\infty$,
\begin{equation}\label{eq3.2}
v(t)=Ae^{-\nu t}+O(e^{-(\nu+\eta)t}),\qquad
v'(t)=-\nu Ae^{-\nu t}+O(e^{-(\nu+\eta)t}).
\end{equation}
The potential profile $W$ satisfies
\begin{equation}\label{eq3.3}
W(t)=
\begin{cases}
\displaystyle\frac{\kappa_NA^p}{a^2-p^2\nu^2}e^{-p\nu t}(1+o(1)),&p\nu<a,\\[1.1ex]
\displaystyle\frac{\kappa_NA^p}{2a}\,te^{-at}(1+o(1)),&p\nu=a,\\[1.1ex]
Be^{-at}(1+o(1)),\quad B>0,&p\nu>a.
\end{cases}
\end{equation}
In each case, $\frac{W'(t)}{W(t)}\to-\min\{a,p\nu\}$.
\end{lemma}

\begin{proof}
By \eqref{eq2.16}, $R=Wv^{p-2}\to0$. For any sufficiently small $\epsilon>0$, comparison in
\[
-v''+(\nu^2-R)v=0
\]
on a half-line gives $v(t)\le C_\epsilon e^{-(\nu-\epsilon)t}$. Indeed, choose $t_0$ so that $R<\nu^2-(\nu-\epsilon)^2$ for $t\ge t_0$. A suitable multiple of $e^{-(\nu-\epsilon)t}$ is then a supersolution that dominates $v$ at $t_0$. Since both functions tend to zero at infinity, the maximum principle gives the estimate.

By $W\le\sqrt{\kappa_N}v$,
\[
0<Wv^{p-1}\le\sqrt{\kappa_N}v^p
\le C_\epsilon e^{-p(\nu-\epsilon)t}.
\]
Choose $\epsilon>0$ so that $p(\nu-\epsilon)>\nu$. The Green representation of $v$ gives \eqref{eq3.2}, with
\[
A=\frac1{2\nu}\int_\R e^{\nu s}W(s)v(s)^{p-1}\,ds>0.
\]
The exponential moment is finite by the preceding estimate and evenness. Differentiating the Green representation gives the remainder for $v'$. In particular, $v^p\in L^1(\R)$.

For $W$, use the identity
\begin{equation}\label{eq3.4}
W(t)=\frac{\kappa_N}{2a}\left[
 e^{-at}\int_{-\infty}^t e^{as}v(s)^p\,ds
+e^{at}\int_t^\infty e^{-as}v(s)^p\,ds\right].
\end{equation}
If $p\nu<a$, both terms have order $e^{-p\nu t}$, and their leading coefficients sum to $\frac{\kappa_NA^p}{a^2-p^2\nu^2}$. If $p\nu=a$, the first integral equals $A^pt+O(1)$, while the second term is $O(e^{-at})$. If $p\nu>a$, the first integral converges, with
\[
B=\frac{\kappa_N}{2a}\int_\R e^{as}v(s)^p\,ds>0.
\]
This proves \eqref{eq3.3}. When we differentiate the representation \eqref{eq3.4}, the endpoint contributions cancel. The same estimates yield the stated limit for $\frac{W'}{W}$.
\end{proof}

Define
\begin{equation}\label{eq3.5}
X=-\frac{v'}v,\qquad Y=-\frac{W'}W,\qquad
R=Wv^{p-2},\qquad d=\nu^2-X^2,\qquad q=\frac Rd.
\end{equation}
By Lemma~\ref{Lem2.5}, $0<X<\nu$ and $0<Y<a$ on $(0,\infty)$, so $d>0$.
The quotient $q$ combines the Riccati equation for $X$ with the logarithmic derivative of $R$. Indeed, $X'=R-d$ and
\[
\frac{R'}R=-Y-(p-2)X,\qquad
\frac{d'}d=2X(1-q),
\]
so that
\[
\frac{q'}q=(2q-p)X-Y.
\]
The Hamiltonian identity below gives $2q-p>0$. Together with the choice $z=-\frac{Y}{a}$, this produces the positive couplings in Lemma~\ref{Lem3.3}.

\begin{lemma}\label{Lem3.2}
The variable $q$ satisfies
\begin{equation}\label{eq3.6}
\lim_{t\to\infty}q(t)=
\begin{cases}
p,&p\nu\le a,\\
\frac{1}{2}\left(p+\frac{a}{\nu}\right),&p\nu>a.
\end{cases}
\end{equation}
\end{lemma}

\begin{proof}
Equations \eqref{eq3.2} and \eqref{eq3.3} give
\begin{equation}\label{eq3.7}
R(t)=Ct^j e^{-mt}(1+o(1)),\qquad
m=(p-2)\nu+\min\{a,p\nu\},\qquad C>0.
\end{equation}
Here $j=1$ if $p\nu=a$ and $j=0$ otherwise. Since $a\ge\nu$ and $p>1$, we have $\min\{a,p\nu\}\ge\nu$, and hence $m\ge(p-1)\nu>0$.

Put $\varepsilon=\nu-X>0$. The Riccati equation gives
\[
\varepsilon'=(\nu+X)\varepsilon-R,
\qquad
\varepsilon(t)=\int_t^\infty
\exp\left(-\int_t^s(\nu+X(\tau))\,d\tau\right)R(s)\,ds.
\]
The boundary term vanishes by \eqref{eq3.2}. Dividing by $R(t)$ and writing $s=t+z$, the integrand converges to $e^{-(2\nu+m)z}$ for each $z\ge0$. By \eqref{eq3.7}, there is a constant $C_1>0$ such that $\frac{R(t+z)}{R(t)}\le C_1(1+z)e^{-\frac{mz}{2}}$ for all sufficiently large $t$ and all $z\ge0$. The other factor is at most one. By dominated convergence,
\begin{equation}\label{eq3.8}
\frac{\varepsilon(t)}{R(t)}\longrightarrow\frac1{2\nu+m},
\qquad
\frac{d(t)}{R(t)}\longrightarrow\frac{2\nu}{2\nu+m}.
\end{equation}
Taking the reciprocal yields \eqref{eq3.6}, including the resonant case.
\end{proof}

Introduce $P=\frac{\kappa_Nv^p}{W}$ and $E=a^2-Y^2$. Dividing \eqref{eq3.1} by $v^2$ yields
\begin{equation}\label{eq3.9}
R(2P-E)=pPd,\qquad P=\frac{qE}{2q-p},\qquad q>\frac p2.
\end{equation}
The last inequality follows from positivity of $q$, $P$, and $E$ in the identity $(2q-p)P=qE$.

\begin{lemma}\label{Lem3.3}
Set $s=\nu t$, $\gamma=\frac{a}{\nu}$, $x(s)=\frac{1}{\nu}X\!\left(\frac{s}{\nu}\right)$, and $z(s)=-\frac{1}{a}Y\!\left(\frac{s}{\nu}\right)$. Also write $q(s)$ for $q\!\left(\frac{s}{\nu}\right)$. With dots denoting derivatives in $s$, the centered profile determines a trajectory in
\[
0<x<1,\qquad -1<z<0,\qquad q>\frac{p}{2}\quad(s>0)
\]
of the system
\begin{equation}\label{eq3.10}
\begin{split}
\dot x&=(1-x^2)(q-1),\\
\dot z&=\gamma(1-z^2)\frac{q-p}{2q-p},\\
\dot q&=q\bigl((2q-p)x+\gamma z\bigr).
\end{split}
\end{equation}
At the center, $x=z=0$. With $q_0=q(0)$, we have
\begin{equation}\label{eq3.11}
\max\left\{1,\frac{p}{2}\right\}<q_0<p,\qquad
v(0)^{2p-2}=\frac{a^2\nu^2q_0^2}{\kappa_N(2q_0-p)},\qquad
W(0)=\nu^2q_0v(0)^{2-p}.
\end{equation}
The system is cooperative and irreducible in this region.
\end{lemma}

\begin{proof}
The Riccati equations are $X'=R-d$ and $Y'=P-E$, while $\frac{R'}{R}=-Y-(p-2)X$. Combining them with \eqref{eq3.9} gives $\frac{q'}{q}=(2q-p)X-Y$. The stated rescaling yields \eqref{eq3.10}.

Since both profiles are even, $d(0)=\nu^2$ and $E(0)=a^2$. Hence $R(0)=q_0\nu^2$ and $P(0)=\frac{q_0a^2}{2q_0-p}$. Since $RP=\kappa_Nv^{2p-2}$, these identities give the initial values in \eqref{eq3.11}. The lower bound $q_0>\frac{p}{2}$ follows from \eqref{eq3.9}. Also $v''(0)=\nu^2v(0)(1-q_0)<0$, so $q_0>1$. Finally,
\[
W(0)=\frac{\kappa_N}{2a}\int_\R e^{-a|s|}v(s)^p\,ds
<\frac{\kappa_N}{a^2}v(0)^p,
\]
which gives $P(0)>a^2$ and therefore $q_0<p$.

The nonzero off-diagonal derivatives of the vector field in \eqref{eq3.10} are
\[
\partial_q\dot x=1-x^2,\qquad
\partial_q\dot z=\frac{\gamma p(1-z^2)}{(2q-p)^2},\qquad
\partial_x\dot q=q(2q-p),\qquad
\partial_z\dot q=\gamma q.
\]
All four derivatives are positive, and each of $x$ and $z$ is coupled in both directions to $q$. This proves cooperativity and irreducibility. These signs also hold at the initial points and along segments joining trajectories in the stated region.
\end{proof}

\begin{proposition}\label{Prop3.4}
For each $0<\nu\le a$, equation \eqref{eq2.12} has exactly one positive $H^1(\R)$ solution up to translation. Its centered representative is the normalized minimizer of \eqref{eq2.10}.
\end{proposition}

\begin{proof}
Existence follows from Proposition~\ref{Prop2.4}, and Lemma~\ref{Lem2.5} centers every positive solution. Suppose two centered solutions have $q_{0,2}>q_{0,1}$. Both trajectories start at $x=z=0$. On each finite interval, their difference solves a linear system with coefficient matrix equal to the average Jacobian along the segment joining them. Lemma~\ref{Lem3.3} gives nonnegative off-diagonal entries and positive couplings. By the comparison principle, the order is preserved and all three inequalities are strict for $s>0$:
\[
x_2>x_1,\qquad z_2>z_1,\qquad q_2>q_1\quad(s>0).
\]
The last equation in \eqref{eq3.10} gives
\[
\frac d{ds}\log\frac{q_2}{q_1}
=2x_2(q_2-q_1)+(2q_1-p)(x_2-x_1)+\gamma(z_2-z_1)>0.
\]
It follows that $\frac{q_2(s)}{q_1(s)}\ge \frac{q_{0,2}}{q_{0,1}}>1$. However, Lemma~\ref{Lem3.2} gives the same positive limit for both $q$ components, so the ratio tends to one. This contradiction shows that the initial parameters coincide.

Equation \eqref{eq3.11} now gives identical initial values, and evenness gives zero initial derivatives. Since the vector field of \eqref{eq2.14} is smooth for $v>0$, uniqueness for the initial value problem proves equality of the profiles.
\end{proof}

\begin{corollary}\label{Cor3.5}
Let $U_\mu(r)=r^{-a}v_\nu(\log r)$, where $v_\nu$ is the centered normalized solution. Then
\[
U_\mu(r)=r^{2-N}U_\mu\!\left(\frac{1}{r}\right)\qquad(r>0).
\]
There is a constant $A_\mu>0$ such that
\begin{equation}\label{eq3.12}
\lim_{r\to0}r^{a-\nu}U_\mu(r)
=\lim_{r\to\infty}r^{a+\nu}U_\mu(r)=A_\mu.
\end{equation}
Moreover, $\frac{rU_\mu'(r)}{U_\mu(r)}$ tends to $-(a-\nu)$ at zero and to $-(a+\nu)$ at infinity.
\end{corollary}

\begin{proof}
Since $v_\nu$ is even, $v_\nu(\log r)=v_\nu(-\log r)$, which is the Kelvin identity. By Lemma~\ref{Lem3.1} and evenness, both limits in \eqref{eq3.12} equal $A_\mu$. Finally,
\[
\frac{rU_\mu'(r)}{U_\mu(r)}=-a+\frac{v_\nu'(\log r)}{v_\nu(\log r)},
\]
and the derivative asymptotics in that lemma give the remaining assertions.
\end{proof}

\section{Nondegeneracy and variational coercivity}

We first justify linearization in $\D$ and then determine the radial and nonradial parts of the kernel. By the Hardy-Littlewood-Sobolev inequality, the map
\[
\mathcal N:\D\longrightarrow L^{\frac{2N}{N+2}}(\R^N),
\qquad \mathcal N(U)=(K_N*|U|^p)|U|^{p-2}U
\]
is well defined.

\Needspace{10\baselineskip}
\begin{lemma}\label{Lem4.1}
At every positive finite energy solution $U$ of \eqref{eq1.1}, the map $\mathcal N$ is Fr\'echet differentiable, with
\[
D\mathcal N(U)h=(p-1)(K_N*U^p)U^{p-2}h
+p\bigl(K_N*(U^{p-1}h)\bigr)U^{p-1}.
\]
The bilinear form associated with $L_\mu$ is continuous on $\D\times\D$. If $1<p<2$, then, for small $h\in\D$,
\begin{equation}\label{eq4.1}
\|\mathcal N(U+h)-\mathcal N(U)-D\mathcal N(U)h\|_{\frac{2N}{N+2}}
\le C_U\|\nabla h\|_2^p.
\end{equation}
\end{lemma}

\begin{proof}
By \eqref{eq2.1}, the local multiplier $(K_N*U^p)U^{p-2}$ belongs to $L^{\frac{N}{2}}$. Since $U^{p-1}h\in L^{\frac{2N}{N+2}}$, the Hardy-Littlewood-Sobolev inequality controls the nonlocal term. For $\phi,\psi\in\D$,
\[
\begin{split}
\left|\int(K_N*U^p)U^{p-2}\phi\psi\right|
&\le\sqrt{\kappa_N}\|U\|_{2^*}^{p-1}\|\phi\|_{2^*}\|\psi\|_{2^*},\\
\left|\int\bigl(K_N*(U^{p-1}\phi)\bigr)U^{p-1}\psi\right|
&\le C\|U\|_{2^*}^{2p-2}\|\phi\|_{2^*}\|\psi\|_{2^*}.
\end{split}
\]
The same H\"older estimates bound the proposed derivative from $\D$ to $L^{\frac{2N}{N+2}}$. The Sobolev and Hardy inequalities then give continuity of the bilinear form.

For $1<p<2$, put $g(t)=|t|^{p-2}t$. The scalar estimates
\begin{equation}\label{eq4.2}
\begin{split}
\bigl||u+h|^p-u^p-pu^{p-1}h\bigr|&\le C_p|h|^p,\\
u\bigl|g(u+h)-u^{p-1}-(p-1)u^{p-2}h\bigr|&\le C_p|h|^p,\\
|g(u+h)-g(u)|&\le C_p|h|^{p-1}
\end{split}
\qquad(u>0)
\end{equation}
follow from Taylor's theorem on $|h|\le\frac{u}{2}$ and direct estimates using $u\le2|h|$ on the complement.

Write $V=K_N*U^p$ and $F_h=|U+h|^p-U^p-pU^{p-1}h$. The nonlinear remainder is exactly
\[
\begin{split}
&V\bigl[g(U+h)-U^{p-1}-(p-1)U^{p-2}h\bigr]\\
&\quad +(K_N*F_h)g(U+h)
+p\bigl(K_N*(U^{p-1}h)\bigr)\bigl[g(U+h)-g(U)\bigr].
\end{split}
\]
By \eqref{eq2.1} and \eqref{eq4.2}, the first term is bounded pointwise by $C|h|^p$. For the second, the Hardy-Littlewood-Sobolev inequality gives $K_N*F_h\in L^{2^*}$, while $g(U+h)\in L^{\frac{N}{2}}$. The same inequality and the last bound in \eqref{eq4.2} control the third term. More explicitly, denote the three displayed terms, in their order, by $I_1,I_2,I_3$, and set $r=\frac{2N}{N+2}$. The estimates are
\[
\begin{aligned}
\|I_1\|_r
&\le C\|h\|_{2^*}^{p},\\
\|I_2\|_r
&\le C\|h\|_{2^*}^{p}
       \bigl(\|U\|_{2^*}+\|h\|_{2^*}\bigr)^{p-1},\\
\|I_3\|_r
&\le C\|U\|_{2^*}^{p-1}\|h\|_{2^*}^{p}.
\end{aligned}
\]
These bounds hold for all sufficiently small real perturbations $h\in\D$. Since $pr=2^*$, all three terms satisfy \eqref{eq4.1}. The remainder is therefore $o(\|\nabla h\|_2)$. For $p\ge2$, the power estimates give continuous differentiability and the same derivative formula.
\end{proof}

Fix $0<\nu\le a$ and the centered normalized solution $v$ from Proposition~\ref{Prop3.4}. Put $W=k*v^p$, $R=Wv^{p-2}$, and $b=v^{p-1}$. The change of variables \eqref{eq2.5} identifies the kernel in $\D$ with the weak $H^1(\Cyl)$ kernel of
\begin{equation}\label{eq4.3}
\mathcal L_\nu=-\partial_t^2-\Delta_\theta+\nu^2-(p-1)R-pbTb.
\end{equation}
Here $(bTb)\phi=b\,T(b\phi)$. The kernel correspondence follows by change of variables for tests supported away from the origin, then by Lemma~\ref{Lem4.1} and density for all test functions in $H^1(\Cyl)$. Write $\Lrad$ for the radial restriction.

\begin{lemma}\label{Lem4.2}
The operator \eqref{eq4.3} is self-adjoint on $L^2(\Cyl)$ with domain $H^2(\Cyl)$ and has essential spectrum $[\nu^2,\infty)$. Its weak $H^1$ kernel equals its operator kernel. In spherical harmonic degree $\ell$, the essential spectrum is $[\nu^2+\lambda_\ell,\infty)$, where $\lambda_\ell=\ell(\ell+N-2)$. In the even radial subspace it is $[\nu^2,\infty)$, and the operator is Fredholm of index zero from $H^2_{\mathrm{even}}(\R)$ to $L^2_{\mathrm{even}}(\R)$.
\end{lemma}

\begin{proof}
Let $A_0=-\partial_t^2-\Delta_\theta+\nu^2$. Fourier transformation and spherical harmonic expansion make $A_0$ self-adjoint on $H^2$, with inverse bounded from $L^2$ to $H^2$. The coefficients $R$ and $b$ are bounded and tend to zero at both ends by \eqref{eq2.16}. Lemma~\ref{Lem2.3} makes $T$ bounded and self-adjoint. Thus the last two terms of \eqref{eq4.3} are bounded self-adjoint perturbations.

To prove compactness of $RA_0^{-1}$, truncate $R$ to a bounded cylinder and apply the local Rellich theorem. The truncation error tends to zero in operator norm because $R$ vanishes uniformly at both ends. The operator $bTb$ is compact on $L^2$ by the same argument, using the $L^2\to H^2$ bound for $T$ and truncating both multipliers. The perturbations therefore preserve the essential spectrum. Restriction to each spherical harmonic subspace, or to the even radial subspace, gives the stated intervals. Since zero lies below the essential spectrum, the maps from the corresponding $H^2$ domains to $L^2$ are Fredholm of index zero.

If $\phi\in H^1(\Cyl)$ belongs to the weak kernel, the lower order terms lie in $L^2$. Thus $A_0\phi\in L^2$, and the inverse estimate gives $\phi\in H^2$. The same argument applies to each restricted kernel.
\end{proof}

A radial kernel element $\phi$ determines
\[
\Psi=\kappa_Np\left(-\frac{d^2}{dt^2}+a^2\right)^{-1}(v^{p-1}\phi).
\]
The pair belongs to $H^2(\R)^2$ and satisfies
\begin{equation}\label{eq4.4}
-\phi''+\nu^2\phi=(p-1)R\phi+v^{p-1}\Psi,
\qquad
-\Psi''+a^2\Psi=\kappa_Np\,v^{p-1}\phi.
\end{equation}
Conversely, an $H^1(\R)^2$ solution of \eqref{eq4.4} gives a radial kernel element. For fixed parameters, $\delta$ denotes differentiation at $\epsilon=0$ under $(v,W)\mapsto(v+\epsilon\phi,W+\epsilon\Psi)$.

\begin{lemma}\label{Lem4.3}
For every $H^1(\R)^2$ solution of \eqref{eq4.4},
\begin{equation}\label{eq4.5}
\delta X(t)\longrightarrow0,
\qquad \frac{\delta q(t)}{q(t)}\longrightarrow0
\quad\text{as }t\to\infty,
\end{equation}
where $X$ and $q$ are defined in \eqref{eq3.5}.
\end{lemma}

\begin{proof}
We show that the relative variations of $R$ and $d$ have the same limit, which gives the limit of $\frac{\delta q}{q}$. This supplies the condition at infinity for the variational comparison in Proposition~\ref{Prop4.4}.

Set $Z=\left(\phi,\frac{\Psi}{\sqrt{\kappa_Np}}\right)$. This symmetrizes \eqref{eq4.4}:
\[
-Z''+\begin{pmatrix}\nu^2&0\\0&a^2\end{pmatrix}Z
=\begin{pmatrix}(p-1)R&\sqrt{\kappa_Np}\,v^{p-1}\\
\sqrt{\kappa_Np}\,v^{p-1}&0\end{pmatrix}Z.
\]
The matrix on the right decays exponentially by Lemma~\ref{Lem3.1} and \eqref{eq3.7}. Fix $0<\beta<\nu$ and choose $t_0$ so that its operator norm is less than $\frac{\nu^2-\beta^2}{2}$ for $t>t_0$. Choose a smooth cutoff $\chi$ that vanishes on $(-\infty,t_0]$ and equals one on $[t_0+1,\infty)$. Testing with $w_L^2Z$, where $w_L=\chi e^{\beta\min\{t,L\}}$ and $L>t_0+1$, gives
\[
\int|(w_LZ)'|^2+\frac{\nu^2-\beta^2}{2}\int w_L^2|Z|^2
\le C_{t_0,\beta}\|Z\|_{H^1(t_0,t_0+1)}^2.
\]
Insert an outer cutoff to justify the test, then remove it for fixed $L$. The terms containing $\chi'$ are supported in $[t_0,t_0+1]$ and bounded independently of $L$. Letting $L\to\infty$ gives $e^{\beta t}Z\in H^1(t_0+1,\infty)$. Estimates on unit intervals, with uniformly bounded coefficients, then yield $Z,Z'=O(e^{-\beta t})$.

Choose $\beta$ close enough to $\nu$ that
\[
m+\beta>\nu,\qquad (p-1)\nu+\beta>\nu,
\]
where $m$ is defined in \eqref{eq3.7}. Such a choice is possible because $m>0$ and $p>1$; when $p<2$, the second condition requires $\beta>(2-p)\nu$. Both terms on the right of the first equation in \eqref{eq4.4} then decay faster than $e^{-\nu t}$. At resonance, the factor $t$ is absorbed by a smaller exponential rate. The Green representation gives $C\in\R$ and $\eta>0$ such that
\[
\phi(t)=Ce^{-\nu t}+O(e^{-(\nu+\eta)t}),\qquad
\phi'(t)=-\nu Ce^{-\nu t}+O(e^{-(\nu+\eta)t}).
\]
The right side of the second equation has leading term $\kappa_NpA^{p-1}Ce^{-p\nu t}$ and a remainder with faster exponential decay. Applying the Green formula \eqref{eq3.4} with $v^p$ replaced by $pv^{p-1}\phi$, and comparing with \eqref{eq3.3}, gives
\[
\frac{\Psi(t)}{W(t)}\longrightarrow
\begin{cases}
\frac{pC}{A},&p\nu\le a,\\
\frac{d_0}{B},&p\nu>a,
\end{cases}
\]
where, in the case $p\nu>a$,
\[
d_0=\frac{\kappa_Np}{2a}
\int_\R e^{as}v(s)^{p-1}\phi(s)\,ds.
\]
This integral is absolutely convergent by the decay estimates above. At resonance, the coefficients of $te^{-at}$ are $\frac{\kappa_NpA^{p-1}C}{2a}$ for $\Psi$ and $\frac{\kappa_NA^p}{2a}$ for $W$. The remaining terms are $O(e^{-at})$, including when $C=0$.

These expansions show that
\[
\vartheta(t):=\frac{\delta R}{R}=\frac\Psi W+(p-2)\frac\phi v\longrightarrow \vartheta_\infty\in\R,
\qquad
\delta X=-\left(\frac\phi v\right)'\longrightarrow0.
\]
To evaluate $\frac{\delta d}{d}$, we need the leading term of $\delta X$ relative to $R$. Linearizing $X'=X^2-\nu^2+R$ gives
\[
(\delta X)'=2X\delta X+R\vartheta,
\qquad (v^2\delta X)'=v^2R\vartheta.
\]
Since $v^2\delta X\to0$, integration from $t$ to infinity yields
\[
\delta X(t)=-v(t)^{-2}\int_t^\infty v(s)^2R(s)\vartheta(s)\,ds.
\]
After division by $R(t)$ and the change $s=t+z$, \eqref{eq3.2} and \eqref{eq3.7} give a constant $C_1>0$ such that, for all sufficiently large $t$ and $z\ge0$,
\[
\left|\frac{v(t+z)^2R(t+z)\vartheta(t+z)}{v(t)^2R(t)}\right|
\le C_1(1+z)e^{-(2\nu+\frac{m}{2})z}.
\]
The bound remains integrable at resonance, so dominated convergence yields
\[
\delta X(t)=-\frac{\vartheta_\infty}{2\nu+m}R(t)+o(R(t)).
\]
Since $\delta d=-2X\delta X$ and $X\to\nu$, \eqref{eq3.8} gives $\frac{\delta d}{d}\to \vartheta_\infty$. Hence
\[
\frac{\delta q}{q}=\frac{\delta R}{R}-\frac{\delta d}{d}\longrightarrow0.
\]
\end{proof}

\begin{proposition}\label{Prop4.4}
For every $0<\nu\le a$, the even radial kernel of $\mathcal L_\nu$ is trivial, and
\[
\Lrad:H^2_{\mathrm{even}}(\R)\longrightarrow L^2_{\mathrm{even}}(\R)
\]
is an isomorphism.
\end{proposition}

\begin{proof}
Let $(\phi,\Psi)\in H^1(\R)^2$ be an even solution of \eqref{eq4.4}. Differentiation of the Hamiltonian using the linearized equations gives $(\delta H)'=0$. The decay estimates in Lemmas~\ref{Lem3.1} and~\ref{Lem4.3} imply $\delta H=0$. Since $\phi'(0)=\Psi'(0)=0$, the vector $(\phi(0),\Psi(0))$ is tangent to $H(v_0,W_0,0,0)=0$, where $v_0=v(0)$ and $W_0=W(0)$.

This curve is regular, since
\[
\partial_{v_0}H=\nu^2v_0(q_0-1)>0.
\]
Formula \eqref{eq3.11} parametrizes its relevant arc by $q_0$, and
\[
\frac d{dq_0}\log v_0
=\frac{q_0-p}{(p-1)q_0(2q_0-p)}<0.
\]
Thus $\delta q_0$ determines every tangent vector. If $\delta q_0=0$, both initial values vanish. The initial derivatives also vanish by evenness, so uniqueness for the linear initial value problem gives the zero pair. For a nonzero pair, change its sign to arrange $\delta q_0>0$.

Before imposing $H=0$, the algebraic identity is
\[
R(2P-E)-pPd=\frac{2pP}{v^2}H.
\]
Since $H=\delta H=0$, differentiation gives the linearized identity in \eqref{eq3.9}. Together with the linearized Riccati equations, this shows that $(\delta x,\delta z,\delta q)$ solves the variational system of \eqref{eq3.10} with initial values $(0,0,\delta q_0)$. To write this system in terms of the relative variation, set
\[
\xi=\delta x,\qquad \zeta=\delta z,\qquad h=\frac{\delta q}{q}.
\]
Then
\[
\begin{aligned}
\dot\xi
&=-2x(q-1)\xi+q(1-x^2)h,\\
\dot\zeta
&=-\frac{2\gamma z(q-p)}{2q-p}\,\zeta
+\frac{\gamma pq(1-z^2)}{(2q-p)^2}\,h,\\
\dot h
&=(2q-p)\xi+\gamma\zeta+2xq\,h,
\end{aligned}
\]
with
\[
(\xi,\zeta,h)(0)=\left(0,0,\frac{\delta q_0}{q_0}\right).
\]
The four nonzero off-diagonal coefficients are positive, so $\xi,\zeta,h>0$ for $s>0$. The last equation then gives $\dot h>0$. Thus $\frac{\delta q}{q}\ge\frac{\delta q_0}{q_0}>0$, contrary to \eqref{eq4.5}. Hence the even radial kernel is trivial. Lemma~\ref{Lem4.2} gives Fredholm index zero on the even subspace, so the operator is bijective and has a bounded inverse.
\end{proof}

Let $\mathfrak q^{\mathrm{rad}}$ be the quadratic form of $\Lrad$ on $H^1(\R)$, with the $L^2(\R,dt)$ pairing. The even coefficients preserve the even and odd subspaces.

\begin{lemma}\label{Lem4.5}
For every $0<\nu\le a$,
\[
\ker\Lrad\cap H^1_{\mathrm{odd}}(\R)=\spn\{v'\}.
\]
\end{lemma}

\begin{proof}
Set $\psi=-v'>0$ on $(0,\infty)$. Translation invariance and $v\in H^3$ show that its odd extension is an $H^1(\R)$ kernel element. Odd $H^1$ functions restrict to $H_0^1(0,\infty)$. On $(0,\infty)$, convolution by $k$ on odd functions is represented by the integral kernel
\[
J(t,s)=k(t-s)-k(t+s)>0\qquad(t,s>0).
\]
For $f\in C_c^\infty(0,\infty)$, let $f_{\mathrm{odd}}$ denote its odd extension to $\R$. Write $f=\psi g$ and test the equation for $\psi$ with $\frac{f^2}{\psi}$. Integration by parts and symmetrization give
\begin{equation}\label{eq4.6}
\begin{split}
\frac12\mathfrak q^{\mathrm{rad}}[f_{\mathrm{odd}}]
={}&\int_0^\infty\psi^2|g'|^2\\
&+\frac p2\iint_{(0,\infty)^2}J(t,s)v(t)^{p-1}v(s)^{p-1}
\psi(t)\psi(s)|g(t)-g(s)|^2\,dt\,ds.
\end{split}
\end{equation}
The test has compact support away from zero. All terms are integrable because $v^{p-1}\psi\in L^2\cap L^\infty$ and $J(t,s)\le k(t-s)+k(t+s)$.

For $f\in H_0^1(0,\infty)$, choose $f_j\in C_c^\infty(0,\infty)$ converging to $f$ in $H^1$. The quadratic form is continuous in this norm. On each compact subinterval, $\psi$ is bounded away from zero, so $\frac{f_j}{\psi}\to\frac{f}{\psi}$ in $H^1$. Pass to a subsequence and apply Fatou's lemma to the nonnegative double integrand. Local convergence followed by exhaustion handles the derivative term. Thus \eqref{eq4.6} extends as a lower bound for $\frac{1}{2}\mathfrak q^{\mathrm{rad}}[f_{\mathrm{odd}}]$. For a kernel element the form vanishes, so $\left(\frac{f}{\psi}\right)'=0$ locally and $f=c\psi$ on $(0,\infty)$. Its odd extension is a multiple of $v'$. By translation invariance, $v_\nu'$ lies in the kernel, which gives the reverse inclusion.
\end{proof}

The one-dimensional arguments above also apply when $a,\nu,\kappa>0$ and $p>1$ are independent parameters with $\nu\le a$. For the system
\[
-v''+\nu^2v=Wv^{p-1},\qquad
-W''+a^2W=\kappa v^p,
\]
the kernel $k(t)=\frac{\kappa}{2a}e^{-a|t|}$ has the properties used in the one-dimensional arguments of Proposition~\ref{Prop2.4} and Lemma~\ref{Lem2.5}. Those arguments give a positive $H^1(\R)^2$ solution. The proofs of Propositions~\ref{Prop3.4} and~\ref{Prop4.4} and Lemma~\ref{Lem4.5} then give uniqueness up to common translation and identify the $H^1(\R)^2$ kernel of the linearized system as $\spn\{(v',W')\}$. The condition $\nu\le a$ gives $W\le\sqrt\kappa\,v$. The critical relation $a(p-1)=2$ is used in the nonradial analysis below.

For the nonradial components, let $\lambda_\ell=\ell(\ell+N-2)$, $\ell=0,1,\ldots$, be the eigenvalues of $-\Delta_\theta$ \cite{AtkinsonHan}. Choose a spherical harmonic of degree $\ell$ with unit $L^2(\SN)$ norm, and let $f(t)$ be its coefficient. The corresponding auxiliary coefficient is
\[
g=\kappa_Np\left(-\frac{d^2}{dt^2}+a^2+\lambda_\ell\right)^{-1}(v^{p-1}f).
\]
For independent $f,g\in H^1(\R)$, define
\begin{equation}\label{eq4.7}
\begin{split}
Q_\ell[f,g]={}&\int_\R\left(f'^2+[\nu^2+\lambda_\ell-(p-1)R]f^2\right)\\
&+\frac1{\kappa_Np}\int_\R\left(g'^2+(a^2+\lambda_\ell)g^2\right)
-2\int_\R v^{p-1}fg.
\end{split}
\end{equation}
For fixed $f$, let $g_f$ be the auxiliary coefficient defined above. Its equation and integration by parts give
\[
Q_\ell[f,g]=Q_\ell[f,g_f]
+\frac1{\kappa_Np}\int_\R
\bigl(|(g-g_f)'|^2+(a^2+\lambda_\ell)|g-g_f|^2\bigr).
\]
Let $\mathfrak q_\ell$ be the quadratic form of the scalar linearization in degree $\ell$. Substitution of $g_f$ gives $\mathfrak q_\ell[f]=Q_\ell[f,g_f]$, and therefore
\begin{equation}\label{eq4.8}
\mathfrak q_\ell[f]=\inf_{g\in H^1(\R)}Q_\ell[f,g].
\end{equation}

\begin{lemma}\label{Lem4.6}
For $0<\nu\le a$ and $\mu=a^2-\nu^2$,
\begin{equation}\label{eq4.9}
\mathfrak q_\ell[f]\ge
\left(\frac{2\mu}{a+\nu}+\lambda_\ell-\lambda_1\right)\|f\|_2^2
\qquad(\ell\ge1).
\end{equation}
In particular, the kernel in $H^1(\Cyl)$ is trivial in every spherical harmonic degree $\ell\ge1$ if $\mu>0$.
\end{lemma}

\begin{proof}
Consider the positive functions
\[
\psi_1=e^{-t}(av-v'),\qquad\chi_1=e^{-t}(aW-W').
\]
Their positivity follows from \eqref{eq2.15}. These functions are the radial coefficients of the translation directions in spherical harmonic degree one. Indeed, $-r^a\partial_{x_j}U(r\theta)=\psi_1(\log r)\theta_j$, with the analogous identity for the Newton potential $V$ and $\chi_1$. Differentiating the Hardy potential produces the nonnegative forcing in the first equation below. Using $a(p-1)=2$ and $\lambda_1=2a+1$, direct differentiation of \eqref{eq2.14} gives
\begin{equation}\label{eq4.10}
\begin{split}
[-\partial_t^2+\nu^2+\lambda_1-(p-1)R]\psi_1-v^{p-1}\chi_1
&=2\mu e^{-t}v,\\
[-\partial_t^2+a^2+\lambda_1]\chi_1-\kappa_Np\,v^{p-1}\psi_1&=0.
\end{split}
\end{equation}
For compactly supported smooth $f,g$, test \eqref{eq4.10} against $\frac{f^2}{\psi_1}$ and $\frac{g^2}{\kappa_Np\chi_1}$, respectively. Integration by parts gives
\begin{equation}\label{eq4.11}
\begin{split}
Q_1[f,g]={}&\int\psi_1^2\left|\left(\frac f\psi_1\right)'\right|^2
+\frac1{\kappa_Np}\int\chi_1^2\left|\left(\frac g\chi_1\right)'\right|^2\\
&+\int v^{p-1}\psi_1\chi_1\left|\frac f\psi_1-\frac g\chi_1\right|^2
+\int\frac{2\mu v}{av-v'}f^2.
\end{split}
\end{equation}
Since the first three terms are nonnegative and $av-v'\le(a+\nu)v$, equation \eqref{eq4.11} gives $Q_1[f,g]\ge\frac{2\mu}{a+\nu}\|f\|_2^2$. Compact support and positivity of $\psi_1,\chi_1$ justify the tests. Continuity of \eqref{eq4.7} and density extend the lower bound to all pairs in $H^1(\R)^2$.

For $\ell\ge1$,
\[
Q_\ell[f,g]=Q_1[f,g]+(\lambda_\ell-\lambda_1)
\left(\|f\|_2^2+\frac{\|g\|_2^2}{\kappa_Np}\right).
\]
Discard the nonnegative $g$ term and use \eqref{eq4.8}. This proves \eqref{eq4.9}.
\end{proof}

\begin{proof}[Proof of Theorem~\ref{Thm1.1}]
Proposition~\ref{Prop2.4} gives a positive normalized extremal. Any positive finite energy solution is radial by Proposition~\ref{Prop2.2}. Its cylinder representative is a positive $H^1$ solution of \eqref{eq2.12}, so Lemma~\ref{Lem2.5} and Proposition~\ref{Prop3.4} identify it with a translate of the centered extremal. Under \eqref{eq2.5}, translation by $\log\lambda$ becomes $U(x)\mapsto\lambda^aU(\lambda x)$.

If $u$ is a nonzero real extremal, then $|u|$ is also an extremal. After amplitude normalization, it is a nonnegative weak solution. Lemma~\ref{Lem2.1} makes it strictly positive and smooth away from the origin. On every compact annulus, $\frac{u}{|u|}\in H^1$ takes values in $\{-1,1\}$ and has zero weak gradient. Connectedness of $\R^N\setminus\{0\}$ implies that $u$ has constant sign. Undoing the amplitude normalization gives the stated family of extremals. Conversely, multiplication by a nonzero constant and dilation preserve the quotient, so every member of this family is an extremal.

By Lemma~\ref{Lem4.2}, every weak kernel element on the cylinder belongs to $H^2(\Cyl)$. Spherical harmonic projections commute with $\mathcal L_\nu$. Lemma~\ref{Lem4.6} rules out nonzero kernel components of spherical harmonic degree $\ell\ge1$ when $\mu>0$, while Proposition~\ref{Prop4.4} and Lemma~\ref{Lem4.5} identify the radial kernel as $\spn\{v'\}$. Differentiating $U(r)=r^{-a}v(\log r)$ gives $r^{-a}v'(\log r)=aU+rU_r=\Lambda U$, proving \eqref{eq1.4}.
\end{proof}

\begin{proof}[Proof of Theorem~\ref{Thm1.2}]
Fix $0<\mu<\mu_H$ and let $U=U_\mu$ be centered and normalized. The form $\mathfrak q_\mu$ corresponds under \eqref{eq2.5} to the form of $\mathcal L_\nu$. Set
\[
\mathcal T_\mu=\{\phi\in\D:\langle U,\phi\rangle_\mu=0\}.
\]

Since $U$ solves the equation, $D(U)=Q_\mu(U)$ and
\[
\langle U,\phi\rangle_\mu=\int_{\R^N}\mathcal N(U)\phi.
\]
The functional $D$ is continuously differentiable on $\D$. Integrating its derivative along $U+t\phi$ and using Lemma~\ref{Lem4.1} gives
\[
D(U+\phi)=D(U)+2p\langle U,\phi\rangle_\mu
+p\int_{\R^N}D\mathcal N(U)[\phi]\,\phi
+o(\|\phi\|_\mu^2).
\]
Since $S_\mu=D(U)^{\frac{p-1}{p}}$, expanding the scalar power gives
\begin{equation}\label{eq4.12}
\begin{split}
Q_\mu(U+\phi)-S_\mu D(U+\phi)^{\frac{1}{p}}
={}&\mathfrak q_\mu[\phi]
+\frac{2(p-1)}{D(U)}\langle U,\phi\rangle_\mu^2
+o(\|\phi\|_\mu^2).
\end{split}
\end{equation}
The left side is nonnegative. Applying this expansion to $t\phi$, with $\phi\in\mathcal T_\mu$, proves $\mathfrak q_\mu[\phi]\ge0$.

Since $\mathcal T_\mu$ has codimension one, the Morse index is at most one. On the other hand,
\[
\mathfrak q_\mu[U]=-2(p-1)D(U)<0.
\]
The Morse index is therefore exactly one. Since the essential spectrum starts above zero, Lemma~\ref{Lem4.2} implies that $\mathcal L_\nu$ has one simple negative eigenvalue. Its eigenfunction is radial and even by the nonnegativity of the odd radial form in \eqref{eq4.6} and the positivity in every positive spherical harmonic degree in \eqref{eq4.9}.

Next suppose $\phi\in\mathcal T_\mu$ and $\mathfrak q_\mu[\phi]=0$. Here $\mathfrak q_\mu(\cdot,\cdot)$ denotes the associated symmetric bilinear form. Nonnegativity on $\mathcal T_\mu$ gives $\mathfrak q_\mu(\phi,z)=0$ for every $z\in\mathcal T_\mu$. By homogeneity,
\[
L_\mu U=-2(p-1)\mathcal N(U),
\qquad \mathfrak q_\mu(\phi,U)=-2(p-1)\langle\phi,U\rangle_\mu=0.
\]
As $\mathcal T_\mu$ and $U$ span $\D$, the vector $\phi$ belongs to the kernel. Theorem~\ref{Thm1.1} identifies it as a multiple of $\Lambda U$. Conversely, dilation invariance gives $\langle U,\Lambda U\rangle_\mu=0$, and $\Lambda U$ is a kernel vector.

To prove \eqref{eq1.5} for this fixed $\mu$, suppose there are $\phi_j$ satisfying both orthogonality conditions with $\|\phi_j\|_\mu=1$ and $\mathfrak q_\mu[\phi_j]\to0$. A subsequence converges weakly in $\D$ to a vector $\phi$ in the same subspace. Put
\[
K_\mu[\psi]=\|\psi\|_\mu^2-\mathfrak q_\mu[\psi].
\]
The compactness argument of Lemma~\ref{Lem4.2} shows that $K_\mu[\phi_j]\to K_\mu[\phi]$. For the local term, this follows from local Rellich compactness on the cylinder and the uniform decay of $R$. For the nonlocal term, it follows from compactness of $bTb$ on $L^2$. Consequently $K_\mu[\phi]=1$ and
\[
0\le\mathfrak q_\mu[\phi]=\|\phi\|_\mu^2-1\le0.
\]
The preceding characterization and orthogonality to $\Lambda U$ imply $\phi=0$, contradicting $K_\mu[\phi]=1$. This proves \eqref{eq1.5}.

For the remainder estimate, use the cylinder isometry and denote its energy inner product by $\langle\cdot,\cdot\rangle_\nu$, with
\[
\|f\|_\nu^2=\int_\Cyl\bigl(|f_t|^2+|\nabla_\theta f|^2+\nu^2f^2\bigr).
\]
Let $f$ represent $u$ and let $v_s(t,\theta)=v_\nu(t+s)$. The set $\mathcal M_\mu$ corresponds to $\{cv_s:c\in\R,\ s\in\R\}$. For $u\ne0$ and $\varrho=d_\mu(u)<\|u\|_\mu$, minimizing first over $c$ gives
\[
\varrho^2=\|f\|_\nu^2-\frac{\displaystyle\sup_{s\in\R}|\langle f,v_s\rangle_\nu|^2}{\|v_\nu\|_\nu^2}.
\]
Approximation with compactly supported functions shows that $s\mapsto\langle f,v_s\rangle_\nu$ is continuous and tends to zero as $|s|\to\infty$. Since $\varrho<\|f\|_\nu$, its absolute value attains a positive maximum. The corresponding minimizing coefficient $c$ is nonzero. Differentiating with respect to $c$ and $s$, using $v_\nu\in H^3(\R)$, gives for $h=f-cv_s$
\[
\langle h,v_s\rangle_\nu=\langle h,v_s'\rangle_\nu=0,
\qquad c^2\|v_\nu\|_\nu^2=\|f\|_\nu^2-\varrho^2.
\]
After translating and dividing by $c$, set $\eta=c^{-1}h(\cdot-s,\cdot)$. Under the cylinder identification, its Euclidean representative satisfies the orthogonality conditions in \eqref{eq1.5}. Moreover,
\[
\|\eta\|_\nu
=\frac{\|v_\nu\|_\nu \varrho}{(\|f\|_\nu^2-\varrho^2)^{\frac{1}{2}}}
\le\frac{\epsilon_\mu}{\sqrt{1-\epsilon_\mu^2}}\|v_\nu\|_\nu
\]
when $\varrho\le\epsilon_\mu\|f\|_\nu$. Choose $\epsilon_\mu$ small enough to apply \eqref{eq4.12} and absorb its remainder into the coercive bound \eqref{eq1.5}. In cylinder variables, the term involving $\langle v_\nu,\eta\rangle_\nu^2$ vanishes by orthogonality. The upper bound follows from continuity of the quadratic form. By translation invariance on the cylinder and homogeneity, the deficit at $f$ equals $c^2$ times the deficit at $v_\nu+\eta$, while $c^2\|\eta\|_\nu^2=\|h\|_\nu^2=\varrho^2$. This proves \eqref{eq1.6} for $u\ne0$, and the case $u=0$ is immediate.
\end{proof}

\begin{lemma}\label{Lem4.7}
If $D(u_n)=1$ and $Q_\mu(u_n)\to S_\mu$, then, along a subsequence,
$r_n^au_n(r_n\,\cdot)$ converges strongly in $\D$ to an element of
$\mathcal M_\mu\setminus\{0\}$ for some $r_n>0$.
\end{lemma}

\begin{proof}
For $F,G\in L^{\frac{2N}{N+2}}(\R^N)$, put
\[
\mathcal D(F,G)=\iint_{\R^N\times\R^N}
\frac{F(x)G(y)}{|x-y|^{N-2}}\,dx\,dy ,
\qquad D(u)=\mathcal D(|u|^p,|u|^p).
\]
Here $\frac{2^*}{p}=\frac{2N}{N+2}>1$.
By the Hardy-Littlewood-Sobolev inequality $\mathcal D$ is continuous on
$L^{\frac{2N}{N+2}}\times L^{\frac{2N}{N+2}}$, and $1=D(u_n)\le C\|u_n\|_{2^*}^{2p}$.
The sequence is bounded in $\D$ by Hardy's inequality. Applying the
concentration lemma of Palatucci and Pisante
(\cite[Section 4]{PalatucciPisante}; see also \cite{LionsLimit,Gerard} for the case of $D^{1,2}$)
with $s=1<\frac{N}{2}$, we obtain $r_n>0$ and
$x_n\in\R^N$ such that, after extraction,
\[
 w_n(y)=r_n^au_n(x_n+r_ny)\rightharpoonup w\ne0 \quad\text{in }\D .
\]
Local Rellich compactness gives almost everywhere convergence after
further extraction.

We record the interaction splitting used below. If $f_n$ is bounded
in $L^{2^*}$ and $f_n\to f$ almost everywhere, the power form of
the Brezis-Lieb lemma \cite[Lemma 2.5]{MorozVS13} gives
\[
 |f_n|^p-|f_n-f|^p\ \longrightarrow\ |f|^p
 \quad\text{in }L^{\frac{2N}{N+2}} .
\]
Here $1<p<2^*$. The sequence $|f_n-f|^p$ is bounded in $L^{\frac{2N}{N+2}}$
and tends to zero almost everywhere, so it is weakly null there.
Expanding the continuous bilinear form $\mathcal D$ therefore gives
\begin{equation}\label{eq4.13}
 D(f_n)=D(f)+D(f_n-f)+o(1).
\end{equation}
In particular, $\sigma:=D(w)\in(0,1]$ since $D(w_n)=1$.

We next exclude $|z_n|\to\infty$, where $z_n=\frac{x_n}{r_n}$ is the
center in units of the concentration scale. Write
\[
 \mathcal H_z(h)=\int_{\R^N}|y+z|^{-2}h(y)^2\,dy,\qquad
 Q_{\mu,z}(h)=\|\nabla h\|_2^2-\mu \mathcal H_z(h).
\]
For fixed $h\in\D$, $\mathcal H_z(h)\to0$ as $|z|\to\infty$. Indeed,
this holds for smooth compactly supported functions by separation of
supports from the pole; approximation in $\D$ gives the result
because $\mathcal H_z(h)^{\frac{1}{2}}\le \mu_H^{-\frac{1}{2}}\|\nabla h\|_2$ uniformly in
$z$. Thus the weighted cross term of $w$ and $w_n-w$ tends
to zero by Cauchy-Schwarz. The Dirichlet cross term tends to zero
by weak convergence. Consequently, if $|z_n|\to\infty$, then
\[
 Q_\mu(u_n)=Q_{\mu,z_n}(w_n)
 =\|\nabla w\|_2^2+Q_{\mu,z_n}(w_n-w)+o(1).
\]
Translation of the integration variable gives
$Q_{\mu,z_n}(h)\ge S_\mu D(h)^{\frac{1}{p}}$. On the fixed nonzero
profile use the free bound $\|\nabla w\|_2^2\ge S_0 D(w)^{\frac{1}{p}}$.
The strict inequality $S_\mu<S_0$ follows by testing with a centered
free extremal $Z$: subtracting $\mu \mathcal H_0(Z)>0$ strictly lowers
its quotient. Taking limits and using \eqref{eq4.13} now gives
\[
 S_\mu\ge S_0 \sigma^{\frac{1}{p}}+S_\mu(1-\sigma)^{\frac{1}{p}}>S_\mu,
\]
because $\sigma>0$ and $t^{\frac{1}{p}}\ge t$ on $[0,1]$. This is impossible.

The centers $z_n$ are therefore bounded (otherwise one could extract
a subsequence tending to infinity in modulus). Pass to $z_n\to z$.
The rescaled functions $\tilde u_n(y)=r_n^au_n(r_n y)=w_n(y-z_n)$
converge weakly to $\tilde u=w(\cdot-z)\ne0$, by strong continuity of
translations on fixed Dirichlet test functions. Again take a subsequence
converging almost everywhere. Dilation preserves $D$ and $Q_\mu$.
Since the pole is now fixed, weak convergence in the equivalent Hilbert
norm $\|\cdot\|_\mu$ gives
$Q_\mu(\tilde u_n)=Q_\mu(\tilde u)+Q_\mu(\tilde u_n-\tilde u)+o(1)$.
Applying the sharp inequality to both terms and using \eqref{eq4.13},
with $D(\tilde u)=\sigma$, yields
\[
 1\ge \sigma^{\frac{1}{p}}+(1-\sigma)^{\frac{1}{p}}.
\]
Strict concavity excludes $0<\sigma<1$. Thus $D(\tilde u)=1$,
$Q_\mu(\tilde u)\ge S_\mu$, and the preceding decomposition of $Q_\mu$ gives
$Q_\mu(\tilde u)=S_\mu$ and $Q_\mu(\tilde u_n-\tilde u)\to0$.
The classification in Theorem~\ref{Thm1.1} gives $\tilde u\in\mathcal M_\mu$.
\end{proof}

\begin{proof}[Proof of Corollary~\ref{Cor1.3}]
Write $\Theta_\mu(u)=Q_\mu(u)-S_\mu D(u)^{\frac{1}{p}}$.
Both $\Theta_\mu$ and $d_\mu^2$ are homogeneous of degree two.
On the unit energy sphere the deficit is bounded below by a positive
constant on $\{d_\mu(u)>\epsilon_\mu\}$, whenever this set is
nonempty. Indeed, if the infimum there were zero, we could choose
$Q_\mu(u_n)=1$, $d_\mu(u_n)>\epsilon_\mu$
and $\Theta_\mu(u_n)\to0$. Then $c_n=D(u_n)^{-\frac{1}{2p}}\to\sqrt{S_\mu}$,
and $c_nu_n$ satisfies the hypotheses of Lemma~\ref{Lem4.7}. Since dilations act
isometrically and preserve $\mathcal M_\mu$, $d_\mu(c_nu_n)\to0$, a contradiction.
As $0\in\mathcal M_\mu$, $d_\mu(u)\le1$ on this sphere. Its positive deficit
away from $\mathcal M_\mu$ and the local lower bound \eqref{eq1.6} of Theorem~\ref{Thm1.2}
near $\mathcal M_\mu$ yield the lower bound in \eqref{eq1.7} by homogeneity.

For the upper bound fix $V=U_{\mu,\lambda}$. Minimization over
the real amplitude gives
\[
 d_\mu(u)^2=Q_\mu(u)-\sup_{\lambda>0}
 \frac{|\langle u,U_{\mu,\lambda}\rangle_\mu|^2}{Q_\mu(U_\mu)} .
\]
The bilinear form $\mathcal D$ is positive semidefinite, hence satisfies
Cauchy-Schwarz: on smooth compactly supported inputs this follows
from $-\Delta |x|^{2-N}=(N-2)|S^{N-1}|\delta_0$ and integration
by parts, and continuity extends it to $L^{\frac{2N}{N+2}}$. The weak equation
for $V$ and H\"older's inequality for the measure
$V(y)^p|x-y|^{2-N}\,dx\,dy$ imply
\[
 |\langle u,V\rangle_\mu|
 \le D(V)^{\frac{p-1}{p}}\mathcal D(V^p,|u|^p)^{\frac{1}{p}}
 \le D(V)^{1-\frac{1}{2p}}D(u)^{\frac{1}{2p}} .
\]
Since $Q_\mu(V)=D(V)$ and $D(V)^{1-\frac{1}{p}}=S_\mu$, the supremum in
the distance formula is at most $S_\mu D(u)^{\frac{1}{p}}$. This proves the
upper bound, including $u=0$.
\end{proof}

The constant $\hat c_\mu$ in Corollary~\ref{Cor1.3} cannot be chosen uniformly as $\mu\to0^+$.
Indeed, for a fixed free bubble $Z_\xi$ centered at $\xi\ne0$,
the deficit tends to zero since
$\bigl(1-\frac{\mu}{\mu_H}\bigr)S_0\le S_\mu\le S_0$. Yet its Dirichlet distance
to the closed radial subspace is positive. The Hardy norm equivalence
then bounds $d_\mu(Z_\xi)$ away from zero for small $\mu$.

The following compactness statement will be used to differentiate the centered family.

\begin{lemma}\label{Lem4.8}
Suppose $\nu_j\to\nu_*>0$ with $0<\nu_j,\nu_*\le a$. Centered minimizers of \eqref{eq2.10} with $B=1$ have a subsequence converging strongly in $H^1(\R)$ to a centered minimizer at $\nu_*$. Moreover, $s_{\nu_j}\to s_{\nu_*}$, and the corresponding normalized solutions converge strongly in $H^2(\R)$ along this subsequence.
\end{lemma}

\begin{proof}
A fixed smooth competitor bounds $s_{\nu_j}$ from above. Since $\nu_j$ is bounded away from zero, the minimizers are uniformly bounded in $H^1$. By \eqref{eq2.13} and local compactness, a subsequence converges weakly in $H^1$ and strongly in $L^{2p}$ to a centered function with $B=1$. If $z_j$ denotes this sequence, then
\[
A_{\nu_j}(z_j)-A_{\nu_*}(z_j)=(\nu_j^2-\nu_*^2)\|z_j\|_2^2\longrightarrow0.
\]
By lower semicontinuity, $s_{\nu_*}\le\liminf_j s_{\nu_j}$. A fixed minimizer at $\nu_*$ gives $\limsup_j s_{\nu_j}\le s_{\nu_*}$. Thus $s_{\nu_j}\to s_{\nu_*}$, the weak limit is a minimizer, and the $A_{\nu_*}$ values converge. The convergence is therefore strong in $H^1$.

The normalization factors $s_{\nu_j}^{\frac{1}{2p-2}}$ converge. For the normalized profiles, $v_j^p\to v_*^p$ in $L^2$ and $v_j^{p-1}\to v_*^{p-1}$ uniformly. Young's inequality gives $k*v_j^p\to k*v_*^p$ in $L^2$. The right sides of \eqref{eq2.12} converge in $L^2$, so the second derivatives converge in $L^2$.
\end{proof}

\begin{proof}[Proof of Corollary~\ref{Cor1.5}]
Lemma~\ref{Lem4.8} and uniqueness in Proposition~\ref{Prop3.4} give continuous dependence of the centered normalized profile in $H^2(\R)$. In this proof, $v_\mu$ denotes $v_\nu$ at $\nu=\sqrt{a^2-\mu}$.

We first prove that $\mu\mapsto v_\mu\in H^2_{\mathrm{even}}(\R)$ is continuously differentiable on $(0,\mu_H)$, with
\begin{equation}\label{eq4.14}
\partial_\mu v_\mu=(\mathcal L^{\mathrm{rad}}_\mu)^{-1}v_\mu,
\end{equation}
where $\mathcal L^{\mathrm{rad}}_\mu$ is the even radial linearization $\Lrad$ at $\nu=\sqrt{a^2-\mu}$.

For the one-dimensional nonlinearity, write $\mathcal N_1(z)=(k*|z|^p)g(z)$ on $\R$, where $g(t)=|t|^{p-2}t$. We differentiate the constructed family through difference quotients. For $1<p<2$, we use the remainder estimate \eqref{eq4.15} at $v_\mu$. At $v=v_\mu$, the decomposition in Lemma~\ref{Lem4.1}, Young's inequality, and $H^2(\R)\subset L^\infty(\R)$ give
\begin{equation}\label{eq4.15}
\|\mathcal N_1(v+h)-\mathcal N_1(v)-D\mathcal N_1(v)h\|_2
\le C_v\|h\|_{H^2}^p\qquad(1<p<2,\ \|h\|_{H^2}\le1).
\end{equation}
Indeed, the term containing $k*v^p$ is bounded in $L^2$ by $C\||h|^p\|_2$, using $k*v^p\le\sqrt{\kappa_N}v$ and \eqref{eq4.2}. The convolution of the $p$th power remainder has $L^2$ norm at most $C\|h\|_{2p}^p$ and multiplies the bounded factor $g(v+h)$. The remaining term satisfies
\[
C\|k*(v^{p-1}h)\|_2\|g(v+h)-g(v)\|_\infty
\le C_v\|h\|_2\|h\|_\infty^{p-1}.
\]
These estimates prove \eqref{eq4.15}, whose right side is $o(\|h\|_{H^2})$. For $p\ge2$, differentiability of the power maps also gives a remainder of order $o(\|h\|_{H^2})$ in $L^2$.

Let $h_\delta=v_{\mu+\delta}-v_\mu$. Subtracting the two equations gives
\[
\mathcal L^{\mathrm{rad}}_\mu h_\delta
=\delta v_{\mu+\delta}+r_\delta,
\qquad \|r_\delta\|_2=o(\|h_\delta\|_{H^2}).
\]
Since the profiles depend continuously on $\mu$, $h_\delta\to0$ in $H^2$. Lemma~\ref{Lem4.2} and Proposition~\ref{Prop4.4} give a bounded inverse on the even subspace. Applying it and absorbing the remainder yields $\|h_\delta\|_{H^2}=O(|\delta|)$. Thus $\frac{r_\delta}{\delta}\to0$ in $L^2$ for increments of either sign within $(0,\mu_H)$. Division by $\delta$ proves \eqref{eq4.14}.

To prove continuity of the derivative, let $\mu_j\to\mu$. Then $b_j=v_{\mu_j}^{p-1}\to b=v_\mu^{p-1}$ uniformly. Strong $L^{2p}$ convergence and $k\in L^2$ give uniform convergence of the convolutions. Since the limit profile is positive, $R_j=(k*v_{\mu_j}^p)v_{\mu_j}^{p-2}$ converges uniformly to $R$ on compact intervals. The bound
\[
0<R_j\le\sqrt{\kappa_N}v_{\mu_j}^{p-1},
\]
together with uniform convergence of $v_{\mu_j}$ and decay of $v_\mu$, controls the coefficients outside a large interval. Hence $R_j\to R$ uniformly on $\R$. Put $L_j=\mathcal L^{\mathrm{rad}}_{\mu_j}$ and $L=\mathcal L^{\mathrm{rad}}_\mu$. Their difference is bounded on $L^2$, with
\[
\begin{split}
\|L_j-L\|_{L^2\to L^2}
\le{}&|\mu_j-\mu|+(p-1)\|R_j-R\|_\infty\\
&+pM_N(\|b_j\|_\infty+\|b\|_\infty)\|b_j-b\|_\infty
\longrightarrow0.
\end{split}
\]
The Neumann series for $I+(L_j-L)L^{-1}$ gives a uniform bound for $L_j^{-1}:L^2_{\mathrm{even}}\to H^2_{\mathrm{even}}$. The identity
\[
L_j^{-1}-L^{-1}=L_j^{-1}(L-L_j)L^{-1}
\]
gives convergence of the inverses in this operator norm. Equation \eqref{eq4.14} proves continuity of the derivative.

The fixed cylinder transformation now gives a $C^1$ family $U_\mu$ in $\D$. At an extremal, the derivative of $\frac{Q_\mu(U)}{D(U)^{\frac{1}{p}}}$ with respect to $U$ vanishes. Differentiation along $U_\mu$ therefore gives \eqref{eq1.11}. The Hardy and Hardy-Littlewood-Sobolev inequalities ensure continuity of the numerator and denominator, and $D(U_\mu)>0$. Hence $S_\mu\in C^1((0,\mu_H))$.
\end{proof}

\section{The best constant and its limiting behavior}

We first prove the comparison bounds, then study the rescaled profiles at the Hardy threshold and use the Fourier representation of their interaction to obtain the second term in the expansion of $S_\mu$. The parameter differentiability established in Section~4 completes the proof of strict concavity.

\begin{proposition}\label{Prop5.1}
The monotonicity, concavity, strict bounds \eqref{eq1.8}, and expansion \eqref{eq1.9} hold for every $N\ge3$.
\end{proposition}

\begin{proof}
For fixed $u\ne0$, the quotient in \eqref{eq1.2} is affine in $\mu$, so $S_\mu$ is concave. If $0\le\mu_1<\mu_2<\mu_H$, use an extremal $U_{\mu_1}$ as a competitor at $\mu_2$, taking a centered standard bubble when $\mu_1=0$. Then
\[
S_{\mu_2}\le S_{\mu_1}-(\mu_2-\mu_1)
\frac{\int|x|^{-2}U_{\mu_1}^2}{D(U_{\mu_1})^{\frac{1}{p}}}
<S_{\mu_1}.
\]

Let $H_N$ be the sharp symmetric Hardy-Littlewood-Sobolev constant for $K_N$; as in the introduction, $S$ denotes the sharp Sobolev constant. Thus $D(u)\le H_N\|u\|_{2^*}^{2p}$. A standard bubble attains equality in both inequalities, so $S_0=H_N^{-\frac{1}{p}}S$. Rearrangement and \eqref{eq2.6} reduce the local Hardy-Sobolev quotient to
\[
\omega_N^{\frac{2}{N}}\frac{A_\nu(v)}{\left(\int|v|^{2^*}\right)^{\frac{2}{2^*}}}.
\]
For $v(t)=w\!\left(\frac{\nu}{a}t\right)$, direct substitution gives
\[
A_\nu(v)=\frac\nu a A_a(w),\qquad
\left(\int|v|^{2^*}\right)^{\frac{2}{2^*}}
=\left(\frac a\nu\right)^{\frac{N-2}{N}}
\left(\int|w|^{2^*}\right)^{\frac{2}{2^*}}.
\]
Taking the infimum therefore gives
\begin{equation}\label{eq5.1}
S^{\mathrm{loc}}_\mu=S\left(\frac\nu a\right)^{\frac{2(N-1)}{N}}
=S\left(1-\frac\mu{\mu_H}\right)^{\frac{N-1}{N}}.
\end{equation}
Combining the sharp inequalities with \eqref{eq5.1} gives the lower bound in \eqref{eq1.8}, obtained in \cite{GuJevnikar} for $N=3,4,5$. Equality would make an extremal for $S_\mu$ a Sobolev bubble, which is bounded at the origin. By Theorem~\ref{Thm1.1}, every nonzero extremal is a constant multiple of a dilation of $U_\mu$, and Corollary~\ref{Cor3.5} shows that it is unbounded at the origin. This is a contradiction.

For the upper bound, use $Z(x)=(1+|x|^2)^{-a}$. Its Hardy and Dirichlet integrals satisfy
\begin{equation}\label{eq5.2}
\frac{\int|x|^{-2}Z^2}{\int|\nabla Z|^2}
=\frac{\mathrm B(a,a)}{(N-2)^2\mathrm B(a+2,a)}
=\frac{4(N-1)}{N(N-2)^2}.
\end{equation}
Here $\mathrm B$ is the beta function. Inserting \eqref{eq5.2} into the quotient gives the upper bound in \eqref{eq1.8}. Equality would make $Z$ an extremal for $S_\mu$, again contradicting Theorem~\ref{Thm1.1} and Corollary~\ref{Cor3.5}. Thus the upper bound is strict. The two bounds have the same linear term at zero and a remainder of order $\mu^2$, which proves \eqref{eq1.9}.
\end{proof}

To study $\mu\uparrow\mu_H$, or equivalently $\nu\downarrow0$, set
\begin{equation}\label{eq5.3}
\rho_\nu(s)=\frac{1}{M_N\nu}k\!\left(\frac{s}{\nu}\right)=\frac{a}{2\nu}\exp\!\left(-\frac{a|s|}{\nu}\right),
\qquad
v(t)=\left(\frac{\nu^2}{M_N}\right)^{\frac{1}{2p-2}}w(\nu t).
\end{equation}
The normalized equation \eqref{eq2.12} becomes
\begin{equation}\label{eq5.4}
-w''+w=(\rho_\nu*w^p)w^{p-1}.
\end{equation}
Write
\[
B_\nu(w)=\int(\rho_\nu*|w|^p)|w|^p,
\qquad
\tau_\nu=\inf_{w\in H^1(\R)\setminus\{0\}}
\frac{\int(w'^2+w^2)}{B_\nu(w)^{\frac{1}{p}}}.
\]
Apart from the amplitude factors, \eqref{eq5.3} multiplies the numerator by $\nu$ and the interaction by $\frac{M_N}{\nu}$. The amplitude factors cancel in the quotient, so
\begin{equation}\label{eq5.5}
s_\nu=M_N^{-\frac{1}{p}}\nu^{1+\frac{1}{p}}\tau_\nu.
\end{equation}
The kernels $\rho_\nu$ have mass one and first absolute moment $\frac{\nu}{a}$. Hence
\begin{equation}\label{eq5.6}
\|\rho_\nu*f\|_2\le\|f\|_2,
\qquad
\|\rho_\nu*f-f\|_2\le\frac\nu a\|f'\|_2
\quad(f\in H^1(\R)).
\end{equation}
The second estimate follows by integrating the translation bound $\|f(\cdot-s)-f\|_2\le|s|\|f'\|_2$ against $\rho_\nu$.

Set
\[
w_*(s)=\frac{p^{\frac{1}{2p-2}}}{\cosh^{\frac{1}{p-1}}((p-1)s)},
\]
and define
\begin{equation}\label{eq5.7}
\begin{split}
I_p&=\frac{p^{\frac{p}{p-1}}\sqrt\pi}{p-1}
\frac{\Gamma\!\left(\frac{p}{p-1}\right)}{\Gamma\!\left(\frac{p}{p-1}+\frac{1}{2}\right)},\\
C_N&=\omega_N^{1-\frac{1}{p}}M_N^{-\frac{1}{p}}I_p^{\frac{p-1}{p}}.
\end{split}
\end{equation}

\begin{proposition}\label{Prop5.2}
Let $v_\nu$ be the centered normalized profile. As $\nu\to0^+$,
\begin{equation}\label{eq5.8}
\left(\frac{M_N}{\nu^2}\right)^{\frac{1}{2p-2}}v_\nu\!\left(\frac{\cdot}{\nu}\right)
\longrightarrow w_*\quad\text{strongly in }H^2(\R).
\end{equation}
Moreover, the two-term expansion \eqref{eq1.10} holds with the constant in \eqref{eq5.7}.
\end{proposition}

\begin{proof}
The limiting variational problem is
\[
\tau_0=\inf_{w\in H^1(\R)\setminus\{0\}}
\frac{\int(w'^2+w^2)}{\left(\int|w|^{2p}\right)^{\frac{1}{p}}}.
\]
The compactness argument of Proposition~\ref{Prop2.4} gives a centered nonnegative minimizer. After normalization, it is a positive solution of $-w''+w=w^{2p-1}$. Multiplication by $w'$ and decay at infinity give $w'^2=w^2-\frac{w^{2p}}{p}$. Thus $w(0)^{2p-2}=p$ and $w'(0)=0$. These initial data determine $w_*$ uniquely. The substitution $y=\tanh((p-1)s)$ yields
\[
\int_\R w_*^{2p}
=\frac{p^{\frac{p}{p-1}}}{p-1}
\mathrm B\left(\frac12,\frac p{p-1}\right)=I_p.
\]
Testing the equation against $w_*$ then gives
\[
\int(w_*'^2+w_*^2)=\int w_*^{2p}=I_p,
\qquad \tau_0=I_p^{\frac{p-1}{p}}.
\]

Young's inequality and the Cauchy-Schwarz inequality give $B_\nu(w)\le\int|w|^{2p}$, hence $\tau_\nu\ge\tau_0$. Testing with $w_*$ and using \eqref{eq5.6} gives $\limsup_{\nu\to0}\tau_\nu\le\tau_0$. Therefore $\tau_\nu\to\tau_0$.

Take any sequence $\nu_j\to0$ and centered minimizers $z_j$ with $B_{\nu_j}(z_j)=1$. Their squared $H^1$ norms equal $\tau_{\nu_j}$ and are uniformly bounded. The Sobolev embedding on $\R$ and the chain rule give a uniform $H^1$ bound for $z_j^p$, since $p>1$. Applying \eqref{eq5.6}, we obtain
\[
\left|\int z_j^{2p}-1\right|
\le\|z_j^p\|_2\|z_j^p-\rho_{\nu_j}*z_j^p\|_2
\le\frac{\nu_j}{a}\|z_j^p\|_2\|(z_j^p)'\|_2\longrightarrow0.
\]
Estimate \eqref{eq2.13} and local compactness give strong $L^{2p}$ convergence along a subsequence to a centered function $z_*$, with $\int z_*^{2p}=1$. Lower semicontinuity and $\tau_{\nu_j}\to\tau_0$ imply
\[
\tau_0\le\int(z_*'^2+z_*^2)\le\lim_j\tau_{\nu_j}=\tau_0.
\]
It follows that $z_*$ minimizes the limiting quotient and that the $H^1$ norms converge. Hence $z_j\to z_*$ strongly in $H^1$.

The normalized profiles $w_j=\tau_{\nu_j}^{\frac{1}{2p-2}}z_j$ converge to $w_*$ in $H^1$. By Proposition~\ref{Prop3.4}, they are the rescaled profiles in \eqref{eq5.8}. To obtain $H^2$ convergence, write
\[
\rho_{\nu_j}*w_j^p-w_*^p
=\rho_{\nu_j}*(w_j^p-w_*^p)+(\rho_{\nu_j}*w_*^p-w_*^p).
\]
The contraction estimate and strong $L^{2p}$ convergence make the first term tend to zero in $L^2$. The second tends to zero by \eqref{eq5.6}. Also $w_j^{p-1}\to w_*^{p-1}$ uniformly, so the right sides of \eqref{eq5.4} converge in $L^2$. The equations give convergence of the second derivatives, hence $w_j\to w_*$ in $H^2$. Every sequence of parameters tending to zero has a subsequence with this limit, so the whole family converges.

In particular, \eqref{eq2.11} and \eqref{eq5.5} give
\[
S_\mu\sim\omega_N^{1-\frac{1}{p}}M_N^{-\frac{1}{p}}I_p^{\frac{p-1}{p}}\nu^{1+\frac{1}{p}}.
\]
Since $\nu^2=\mu_H-\mu$ and $\frac{1}{2}\left(1+\frac{1}{p}\right)=\frac{N}{N+2}$, this proves the leading asymptotic.

To obtain the second term, write $B_0(w)=\int_\R|w|^{2p}$ and $A_1(w)=\int_\R(w'^2+w^2)$, in the notation of \eqref{eq2.10}. With the unitary Fourier transform, convolution by $\rho_\nu$ has multiplier $\frac{a^2}{a^2+\nu^2\xi^2}$. Hence, for $f=|w|^p\in H^1(\R)$, Plancherel's identity gives
\begin{equation}\label{eq5.9}
\frac{B_0(w)-B_\nu(w)}{\nu^2}
=\int_\R\frac{\xi^2}{a^2+\nu^2\xi^2}|\widehat f(\xi)|^2\,d\xi.
\end{equation}
The integral on the right converges to $a^{-2}\|f_*'\|_2^2$ whenever $f_\nu\to f_*$ in $H^1$. Replacing $f_\nu$ by $f_*$ changes the integral by at most
\[
a^{-2}(\|f_\nu'\|_2+\|f_*'\|_2)\|f_\nu'-f_*'\|_2,
\]
and the limit for $f_*$ follows by dominated convergence.

Let $w_\nu$ denote the rescaled normalized solution, which minimizes $\frac{A_1}{B_\nu^{\frac{1}{p}}}$ and has just been shown to converge to $w_*$ in $H^2$. Uniform convergence, strong $L^{2p}$ convergence, and $(w_\nu^p)'=p w_\nu^{p-1}w_\nu'$ imply $w_\nu^p\to w_*^p$ in $H^1$. Consequently,
\[
e_\nu:=\frac{B_0(w_\nu)-B_\nu(w_\nu)}{\nu^2B_0(w_\nu)}
\longrightarrow e_*:=\frac{\|(w_*^p)'\|_2^2}{a^2I_p}.
\]
As $\frac{A_1(w_\nu)}{B_0(w_\nu)^{\frac{1}{p}}}\ge\tau_0$,
\[
\tau_\nu
=\frac{A_1(w_\nu)}{B_0(w_\nu)^{\frac{1}{p}}}(1-\nu^2e_\nu)^{-\frac{1}{p}}
\ge\tau_0\left(1+\frac{e_*}{p}\nu^2+o(\nu^2)\right).
\]
Testing with $w_*$ and using \eqref{eq5.9} gives the matching upper bound. Thus
\begin{equation}\label{eq5.10}
\tau_\nu=\tau_0\left(1+\frac{\nu^2}{pa^2}
\frac{\|(w_*^p)'\|_2^2}{\|w_*^p\|_2^2}+o(\nu^2)\right).
\end{equation}
To compute the coefficient, note that $\frac{w_*'}{w_*}=-\tanh((p-1)s)$. The substitution $y=\tanh((p-1)s)$ yields
\[
\frac{\|(w_*^p)'\|_2^2}{\|w_*^p\|_2^2}
=p^2\frac{\mathrm B\!\left(\frac{3}{2},\frac{p}{p-1}\right)}{\mathrm B\!\left(\frac{1}{2},\frac{p}{p-1}\right)}
=\frac{p^2(p-1)}{3p-1}.
\]
Combining \eqref{eq5.10}, \eqref{eq5.5}, and \eqref{eq2.11}, the relative coefficient is
\[
\frac{p(p-1)}{a^2(3p-1)}
=\frac{8(N+2)}{(N-2)^3(N+4)}=b_N.
\]
Since $\nu^2=\mu_H-\mu$, this proves \eqref{eq1.10}.
\end{proof}

\begin{proof}[Proof of Theorem~\ref{Thm1.4}]
Proposition~\ref{Prop5.1} gives the strict bounds, monotonicity, concavity, and expansion at zero. Proposition~\ref{Prop5.2} gives the expansion at the Hardy threshold, including its explicit second coefficient.

It remains to prove $C^2$ regularity and $S_\mu''<0$ on $(0,\mu_H)$. Let
\[
\widehat U_\mu=D(U_\mu)^{-\frac{1}{2p}}U_\mu,
\qquad z_\mu=\partial_\mu\widehat U_\mu.
\]
We shall prove
\begin{equation}\label{eq5.11}
S_\mu''=-2\mathfrak q_\mu[z_\mu]<0.
\end{equation}

By Corollary~\ref{Cor1.5}, the map $\mu\mapsto\widehat U_\mu\in\D$ is $C^1$, and $D(\widehat U_\mu)=1$. Formula \eqref{eq1.11} becomes
\[
S_\mu'=-\int_{\R^N}|x|^{-2}\widehat U_\mu^2.
\]
The right side is continuously differentiable by the Hardy inequality, so $S_\mu\in C^2$ and
\[
S_\mu''=-2\int_{\R^N}|x|^{-2}\widehat U_\mu z_\mu.
\]
Using $\mathcal N$ for the Euclidean nonlinearity of Section~4, the normalized Euler-Lagrange equation is
\[
(-\Delta-\mu|x|^{-2})\widehat U_\mu=S_\mu\mathcal N(\widehat U_\mu).
\]
If $\alpha_\mu=D(U_\mu)^{-\frac{1}{2p}}$, then $S_\mu\alpha_\mu^{2p-2}=1$. Since $\mathcal N$ is homogeneous of degree $2p-1$, its Fr\'echet differentiability at $U_\mu$ transfers to every positive multiple of $U_\mu$. Lemma~\ref{Lem4.1} and this homogeneity give
\[
S_\mu D\mathcal N(\widehat U_\mu)=D\mathcal N(U_\mu).
\]
Differentiating the equation in the dual space $(\D)^*$ yields
\[
L_\mu z_\mu=|x|^{-2}\widehat U_\mu+S_\mu'\mathcal N(\widehat U_\mu).
\]
The derivative of $D(\widehat U_\mu)=1$ is $\int\mathcal N(\widehat U_\mu)z_\mu=0$. Pairing the last equation with $z_\mu$ proves the equality in \eqref{eq5.11}.

Differentiating the normalization gives $\langle U_\mu,z_\mu\rangle_\mu=0$. The cylinder representative of $z_\mu$ is even and that of $\Lambda U_\mu$ is odd, so $\langle z_\mu,\Lambda U_\mu\rangle_\mu=0$. Estimate \eqref{eq1.5} therefore gives $\mathfrak q_\mu[z_\mu]\ge\gamma_\mu\|z_\mu\|_\mu^2$, and it remains to show $z_\mu\ne0$. Suppose $z_\mu=0$. The differentiated equation would then imply
\[
|x|^{-2}\widehat U_\mu=-S_\mu'\mathcal N(\widehat U_\mu).
\]
Write $v_\mu$ for the centered profile at $\nu=\sqrt{a^2-\mu}$ and set $\widehat v_\mu(t)=e^{at}\widehat U_\mu(e^t)=\alpha_\mu v_\mu(t)$. In cylinder variables, the last identity makes $(k*\widehat v_\mu^p)\widehat v_\mu^{p-2}$ the positive constant $-\frac{1}{S_\mu'}$. This contradicts
\[
(k*\widehat v_\mu^p)\widehat v_\mu^{p-2}=\alpha_\mu^{2p-2}R\longrightarrow0
\]
at both ends, by \eqref{eq2.16}. Thus $z_\mu\ne0$, and the strict inequality in \eqref{eq5.11} follows.

By \eqref{eq1.8}, $S_\mu$ is continuous at zero. For $0<\mu_2<\mu_H$, the difference between $\mu\mapsto S_\mu$ and its chord over $[0,\mu_2]$ is concave, vanishes at the endpoints, and is strictly concave in the interior. It is therefore positive on $(0,\mu_2)$, which proves strict concavity on $[0,\mu_H)$.
\end{proof}

\section*{Acknowledgments}
{\sloppy The authors were supported by the National Natural Science Foundation of China (12301145, 12561020, 12261107) and Yunnan Fundamental Research Projects (202401AU070123, 202601AT070048).\par}

\medskip
{\bf AI Declaration:} The authors used GPT-6 (OpenAI) for literature searches, exploratory mathematical discussions, proof checking, and English language editing. All mathematical statements, proofs, and references were independently checked by the authors, who take full responsibility for the content of the paper.

\medskip
{\bf Data availability:} Data sharing is not applicable to this article as no new data were created or analyzed in this study.

\medskip
{\bf Conflict of interest:} The authors declare that they have no conflict of interest.

\end{document}